\documentclass[12pt]{article}
\usepackage{latexsym,amssymb,amsmath}
\begin{document}
\baselineskip=20pt

\parskip=10pt
\def\dis{\displaystyle}
\newcommand{\rd}{\mbox{Rad}}
\newcommand{\kn}{\mbox{ker}}
\newcommand{\psp}{\vspace{0.4cm}}
\newcommand{\pse}{\vspace{0.2cm}}
\newcommand{\ptl}{\partial}
\newcommand{\dlt}{\delta}
\newcommand{\Dlt}{\Delta}
\newcommand{\sgm}{\sigma}
\newcommand{\al}{\alpha}
\newcommand{\be}{\beta}
\newcommand{\G}{\Gamma}
\newcommand{\gm}{\gamma}
\newcommand{\lmd}{\lambda}
\newcommand{\td}{\tilde}
\newcommand{\ad}{\mbox{ad}}
\newcommand{\stl}{\stackrel}
\newcommand{\ol}{\overline}
\newcommand{\es}{\epsilon}
\newcommand{\la}{\langle}
\newcommand{\ra}{\rangle}
\newcommand{\vf}{\varphi}
\newcommand{\vsi}{\varsigma}
\newcommand{\ves}{\varepsilon}
\newcommand{\vt}{\vartheta}
\newcommand{\wt}{\mbox{wt}\:}
\newcommand{\sym}{\mbox{sym}}
\newcommand{\for}{\mbox{for}}
\newcommand{\mbb}{\mathbb}

\def\der{\mbox{der}}
\def\a{\alpha}
\def\b{\beta}
\def\sF{\hbox{$\sc I\hskip -3.5pt F$}}
\def\Z{\hbox{$Z\hskip -5.2pt Z$}}
\def\sZ{\hbox{$\sc Z\hskip -4.2pt Z$}}
\def\Q{\hbox{$Q\hskip -5pt\vrule height 6pt depth 0pt\hskip 6pt$}}
\def\R{\hbox{$I\hskip -3pt R$}}
\def\C{\hbox{$C\hskip -5pt \vrule height 6pt depth 0pt \hskip 6pt$}}
\def\J {\vec J}
\def\d{\delta}
\def\D{\Delta}
\def\g{\gamma}
\def\G{\Gamma}
\def\l{\lambda}
\def\L{\Lambda}
\def\o{\omiga}
\def\p{\psi}
\def\Si{\Sigma}
\def\si{\sigma}
\def\sc{\scriptstyle}
\def\ssc{\scriptscriptstyle}
\def\dis{\displaystyle}
\def\cl{\centerline}
\def\DD{{\cal D}}
\def\ll{\leftline}
\def\rl{\rightline}
\def\sF{\hbox{$\sc I\hskip -2.5pt F$}}
\def\ol{\overline}
\def\ul{\underline}
\def\wt{\widetilde}
\def\wh{\widehat}
\def\rar{\rightarrow}
\def\Rar{\Rightarrow}
\def\lar{\leftarrow}
\def\Lar{\Leftarrow}
\def\rla{\leftrightarrow}
\def\Rla{\Leftrightarrow}
\def\bs{\backslash}
\def\hs{\hspace*}
\def\vs{\vspace*}
\def\rb{\raisebox}
\def\ra{\rangle}
\def\la{\langle}
\def\Rad{\mbox{Rad}}
\def\SS{\hbox{$S\hskip -6.2pt S$}}
\def\hi{\hangindent}
\def\ha{\hangafter}

\def\Bbb#1{{\mbox{b$\!\!\!\!$}\cal #1}}
\def\AA{{\cal A}}
\def\BB{{\cal B}}
\def\CC{{\cal C}}
\def\JJ{{\cal J}}
\def\vi{\vec i}
\def\vj{\vec j}
\def\vk{\vec k}
\def\vm{\vec m}
\def\ii{{\bf i}}
\def\jj{{\bf j}}
\def\kk{{\bf k}}
\def\sone{{1\hskip -5.5pt 1}}
\def\one{{1\hskip -6.5pt 1}}
\def\sJ{J}
\def\th{\theta}
\def\isom{\mbox{${}^{\cong}_{\rar}$}}

\begin{center}{\Large \bf Structure of Contact Lie Algebras Related to}
\end{center}
\begin{center}{\Large \bf Locally-Finite Derivations}\end{center}

\begin{center}{\large Yucai Su$^{\ast}$ and  Xiaoping Xu$^{\dag}$}\end{center}

* Department of Applied Mathematics, Shanghai Jiaotong University, 1954 Huashan Road, Shanghai 200030, P. R. China.

\dag Department of Mathematics, The Hong Kong University of Science \& Technology, Clear Water Bay, Kowloon, Hong Kong, P. R. China\footnote{Research supported by Hong Kong RGC Competitive Earmarked Research Grant HKUST6133/00p}

\vspace{0.3cm}

\begin{abstract}{Classical contact Lie algebras are the fundamental algebraic structures on the manifolds of contact elements of configuration spaces in classical mechanics. In this paper, we determine the structure of the currently largest known category  of contact simple Lie algebras introduced earlier by the second author. These algebras are in general not finitely-graded.}\end{abstract}

\section{Introduction}

A {\it contact element} to an $n$-dimensional smooth manifold at some point is an $(n-1)$-dimensional hyperplane tangent to the manifold at the point. The set of all contact elements of an $n$-dimensional manifold has a natural smooth manifold structure of dimension $2n-1$. Classical contact Lie algebras are the fundamental algebraic structures on the manifolds of contact elements of configuration spaces in classical mechanics (cf.~[A]). Moreover, some contact Lie algebras are generated by certain  quadratic  conformal algebras (cf.~[X2]). The representations theory of the Lie algebras generated by  conformal algebras can be viewed as the algebraic entity of two-dimensional quantum field theory (cf.~[K3]).

A {\it Poisson algebra} is a vector space ${\cal A}$ with two algebraic operations $\cdot$ and $[\cdot,\cdot]$ such that $({\cal A},\cdot)$ forms a commutative associative algebra, $({\cal A},[\cdot,\cdot])$ forms a Lie algebra and the following compatibility condition holds:
$$[u, v\cdot w]=[u,v]\cdot w+v\cdot[u,w]\qquad\for\;\;u,v,w\in {\cal A}.\eqno(1.1)$$
Denote by $\mbb{Z}$ the ring of integers. A {\it (generalized) contact Lie algebra} is a Lie algebra structure on a commutative associative algebra $( {\cal A},\cdot)$  whose Lie bracket is of the form
$$[u,v]=[u,v]_0+(2-\ptl)(u)\ptl_0(v)-\ptl_0(u)(2-\ptl)(v)\qquad\for\;\; u,v\in{\cal A},\eqno(1.2)$$
where $({\cal A},\cdot,[\cdot,\cdot]_0)$ forms a Poisson algebra, $\ptl$ and $\ptl_0$ are mutually commutative derivations of $({\cal A},\cdot)$ such that
$$\ptl([u,v]_0)=[\ptl(u),v]_0+[u,\ptl(v)]_0-2[u,v]_0,\;\;\ptl_0([u,v]_0)=[\ptl_0(u),v]_0+[u,\ptl_0(v)]_0\eqno(1.3)$$
for $u,v\in{\cal A}.$ For a contact Lie algebra structure on  a commutative associative algebra $( {\cal A},\cdot)$,  we define
$$P_u(v,w)=[u,v\cdot w]-[u,v]\cdot w-v\cdot[u,w]\qquad \for\;\; u,v,w\in {\cal A}.\eqno(1.4)$$
Then the following equation holds:
$$P_{u\cdot v}(w_1,w_2)=u\cdot P_v(w_1,w_2)+v\cdot P_u(w_1,w_2) \qquad \for\;\; u,v,w_1,w_2\in {\cal A}.\eqno(1.5)$$
The supersymmetric version of the above equation is the main axiom of  weak Frobenius manifold (cf.~[HM]). Frobenius manifold was introduced by Dubrovin [D] , in connection with topological field theories. In this paper, we determine the isomorphism classes of the currently largest known category  of contact simple Lie algebras introduced by the second author [X1].

 A Lie algebra ${\cal G}$ is called {\it finitely-graded} if ${\cal G}=\bigoplus_{\al\in\G}{\cal G}_{\al}$ is a $\G$-graded vector space for some abelian group $\G$ such that
$$\dim\:{\cal G}_{\al}<\infty,\;\;[{\cal G}_{\al},{\cal G}_{\be}]\subset{\cal G}_{\al+\be}\qquad\for\;\;\al,\be\in \G.\eqno(1.6)$$
Finitely-graded contact simple Lie algebras  have been studied by Kac [K1], [K2], Osborn [O], and Osborn and Zhao [OZ]. The most important feature of the contact simple Lie algebras constructed in [X1] is that they are in general not finitely-graded. Below, we shall give a more technical introduction.

Throughout this paper, we denote by $\mbb{F}$ a field with characteristic 0. All the vector spaces (algebras) are assumed over $\mbb{F}$. Moreover, we denote
by $\mbb{N}$ the additive semi-group of nonnegative integers. When the context is clear, we shall omit the symbol for associative algebraic operation in a product. For $m,n\in\mbb{N}$, we shall use the following notation of indices
$$\ol{m,n}=\left\{\begin{array}{ll}\{m,m+1,...,n\}&\mbox{if}\;m\le n
\\ \emptyset&\mbox{if}\;m>n.\end{array}\right.
\eqno(1.7)$$

Let ${\cal  A}=\mbb{F}[t_1,t_2,...,t_n]$ be the algebra of polynomials in $n$ variables. Recall that a {\it derivation} $\ptl$  of ${\cal A}$ is a linear transformation of ${\cal A}$ such that
$$\ptl(uv)=\ptl(u)v+u\ptl(v)\qquad\for\;\;u,v\in{\cal A}.\eqno(1.8)$$
The typical derivations are $\{\ptl_{t_1},\ptl_{t_2},...,\ptl_{t_n}\}$,  the operators of taking partial derivatives. The space $\mbox{Der}\:{\cal A}$ of all the derivations of ${\cal A}$ forms a Lie algebra with respect to the commutator. Identifying the elements of
${\cal A}$ with their corresponding multiplication operators, we have
$$\mbox{Der}\:{\cal A}=\sum_{i=1}^n{\cal A}\ptl_{t_i},\eqno(1.9)$$
which forms a simple Lie algebra. The Lie algebra $\mbox{Der}\:{\cal A}$ is called a {\it Witt algebra of rank} $n$, usually denoted as ${\cal W}(n,\mbb{F})$. The Lie algebra ${\cal W}(n,\mbb{F})$ acts on the Grassmann algebra $\hat{\cal A}$ of differential forms on ${\cal A}$ as follows.
$$\ptl(df)=d(\ptl(f)),\;\;\ptl(\omega\wedge \nu)=\ptl(\omega)\wedge\nu +\omega\wedge\ptl(\nu)\eqno(1.10)$$
for $f\in{\cal A},\;\omega,\nu\in\hat{\cal A},\;\ptl\in {\cal W}(n,\mbb{F})$.
Assume that $n=2k+1$ is an odd integer. We let
\begin{eqnarray*}{\cal K}(2k+1,\mbb{F})&=&\{\ptl\in {\cal W}(n,\mbb{F})\mid \ptl(dt_{2k+1}+\sum_{i=1}^k(t_idt_{k+i}-t_{k+i}dt_i)\\& &\in {\cal A}(dt_{2k+1}+\sum_{i=1}^k(t_idt_{k+i}-t_{k+i}dt_i)\}.\hspace{5.2cm}(1.11)\end{eqnarray*}
The subspace  ${\cal K}(2k+1,\mbb{F})$ forms a simple Lie subalgebra , which is called a {\it classical contact Lie algebra}.

To see that the Lie algebra ${\cal K}(2k+1,\mbb{F})$ is isomorphic to a Lie algebra with the Lie bracket of the form (1.2), we set
$$\es(i)=1,\;\;\es(j)=-1,\;\;\ol{i}=k+i,\;\;\ol{j}=j-k\qquad\for\;\;i\in\ol{1,k},\;j\in\ol{k+1,2k}.\eqno(1.12)$$
Moreover, we define
$$D_K(f)=\sum_{i=1}^{2k}(t_i\ptl_{t_n}(f)+\es(\ol{i})\ptl_{t_{\ol{i}}}(f))\ptl_{t_i}+(2f-\sum_{i=1}^{2k}t_{\ol{i}}(\es(\ol{i})t_i\ptl_{t_n}(f)+\ptl_{t_{\ol{i}}}(f)))\ptl_{t_n}\eqno(1.13)$$
for $f\in{\cal A}$. Then the map $D_K: f\rar D_k(f)$ is a linear isomorphism from ${\cal A}$ to ${\cal K}(2k+1,\mbb{F})$ (e.g., cf.~[SF]). Furthermore, we define
$$[f,g]_0=\sum_{i=1}^{2k}\es(i)\ptl_{t_i}(f)\ptl_{t_{\ol{i}}}(g)\qquad\for\;\;f,g\in{\cal A}\eqno(1.14)$$
and
$$\ptl=\sum_{i=1}^{2k}t_i\ptl_{t_i},\;\;\;\;\ptl_0=\ptl_{t_n}.\eqno(1.15)$$
It can be verified that
$$[D_K(f),D_K(g)]=D_K([f,g]_0+(2-\ptl)(f)\ptl_0(g)-\ptl_0(f)(2-\ptl)(g))\eqno(1.16)$$
for $f,g\in{\cal A}$ (e.g., cf.~[SF]). Thus the Lie bracket of ${\cal K}(2k+1,\mbb{F})$ is essentially of the form (1.2).

Define the grading
$$ ({\cal K}(2k+1,\mbb{F}))_m=\mbox{Span}\:\{D_K(t_1^{\ell_1}t_2^{\ell_2}\cdots t_n^{\ell_n})\mid \ell_i\in\mbb{N},\;2\ell_n+\sum_{i=1}^{2k}\ell_i=m+2\}\eqno(1.17)$$
for $-2\leq m\in\mbb{Z}$. Then  ${\cal K}(2k+1,\mbb{F})$ is a finitely-graded Lie algebra.

The contact  Lie algebras constructed in [X1] are as follows. Let ${\cal A}$ be a commutative associative algebra with an identity element and let $\{\ptl_0,\ptl_1,...,\ptl_{2k}\}$ be $2k+1$ mutually commutative derivations of ${\cal A}$, where $k$ is a positive integer.
Pick any elements
$$\{\xi_0,\xi_1,\xi_2,...,\xi_k\}\subset {\cal A}\eqno(1.18)$$
such that
$$\ptl_i(\xi_0)=\ptl_{k+i}(\xi_0)=\ptl_0(\xi_j)=\ptl_i(\xi_j)=\ptl_{k+i}(\xi_j)=0\eqno(1.19)$$
for $i,j\in\ol{1,k},\;i\neq j.$
We define an algebraic operation $[\cdot,\cdot]_0$ on ${\cal A}$ by
$$[u,v]_0=\sum_{i=1}^k\xi_i(\ptl_i(u)\ptl_{k+i}(v)-\ptl_{k+i}(u)\ptl_i(v))\qquad\for\;\;u,v\in{\cal A}.\eqno(1.20)$$
Then $({\cal A},\cdot, [\cdot,\cdot]_0)$ forms a Poisson algebra. Take a derivation $\ptl$ of ${\cal A}$ such that the first equation in (1.3) holds and
$$\ptl\ptl_0=\ptl_0\ptl,\;\;\;\ptl(\xi_0)=0.\eqno(1.21)$$
The Lie bracket of the contact Lie algebra in [X1] is defined as
$$[u,v]=[u,v]_0+\xi_0[(2-\ptl)(u)\ptl_0(v)-\ptl_0(u)(2-\ptl)(v)]\qquad\for\;\; u,v\in{\cal A},\eqno(1.22)$$
which is of the form (1.2) if we replace $\xi_0\ptl_0$ by $\ptl_0$.

A linear transformation $T$ on a vector space $V$ is called {\it locally-finite} if
$$\dim(\mbox{Span}\:\{T^m(u)\mid m\in\mbb{N}\})<\infty\qquad\for\;\;u\in V.\eqno(1.23)$$
In [X1], the second author proved that the Lie algebra ${\cal A}$ with Lie bracket (1.22) is simple when ${\cal A}$ is a certain semigroup algebra, $\xi_i$ for $i\in\ol{0,k}$ are invertible, $\ptl_j$ for $j\in\ol{0,2k}$ are locally-finite, $\ptl$ is diagonalizable and some other distinguishable conditions among $\{\ptl,\ptl_j\mid j\in\ol{0,2k}\}$ hold. We refer to [SXZ] for the classification of derivation-simple algebras when the derivations are locally-finite. The Lie algebra $({\cal A},[\cdot,\cdot])$ is in general not finitely-graded. The contact simple Lie  algebras studied in [K1], [O] and [OZ] are special finitely-graded  cases of those in [X1]

In Section 2, we shall rewrite the presentations of contact simple
Lie algebras given in [X1] up to certain relatively obvious
isomorphisms, which we call {\it normalized forms}. In Section 3,
we shall present five lemmas that will be used in the proof of the
main theorem on the isomorphism classes of the normalized contact
simple Lie algebras. Section 4 is devoted to the presentation of
the main theorem and its proof.

\section{Normalized form of contact Lie algebras}

In this section, we shall present the normalized form of contact simple  Lie algebras introduced in [X1].

Let
$$\vec\ell=(\ell_1,...,\ell_6)\in \mbb{N}\:^6\mbox{ such that }\sum_{p=1}^6\ell_p>0.
\eqno(2.1)$$
Denote
$$
\iota_i=\ell_1+\ell_2+...+\ell_i\qquad \for\ \ i=1,2,...,6,
\eqno(2.2)$$
and define $6$ index sets
$$
I_i=\ol{\iota_{i-1}+1,\iota_i}\qquad \for\ \ i=1,2,...,6,
\eqno(2.3)$$
where we treat $\iota_{-1}=0$. Denote by
$I_{i,j}$ the union from the $i$th index set to the $j$th index:
$$
I_{i,j}=\bigcup_{i\le p\le j}I_p=\ol{\iota_{i-1}+1,\iota_j}\qquad
\for\ 1\le i\le j\le6.
\eqno(2.4)$$
Set
$$I=I_{1,6}=\ol{1,\iota_6},\qquad J=\ol{1,2\iota_6}.\eqno(2.5)$$
Moreover, we denote
$$
\wh K=\{0\}\cup K\qquad\for\;\;K\subset J.\eqno(2.6)$$
Define the map \rb{5pt}{$\ol{\ }$}\,$:J\rar J$ to be the index
shifting of $\iota_6$ steps in $J$, i.e.,
$$
\ol p=\left\{\begin{array}{ll}p+\iota_6&\mbox{if \ }p\in\ol{1,\iota_6},\vs{4pt}\\
p-\iota_6&\mbox{if \ }p\in\ol{\iota_6+1,2\iota_6}.\\
\end{array}\right.\eqno(2.7)$$
We also define $\ol 0=0$ for convenience.
For any subset $K$ of $\wh J$, we denote
$$\ol K=\{\ol p\,|\,p\in K\}.\eqno(2.8)$$
Thus $J=I\cup\ol I$. Set
$$
J_i=I_i\cup\ol I_i,\ \ J_{i,j}=I_{i,j}\cup\ol I_{i,j}.
\eqno(2.9)$$

An  element of the vector space $\mbb{F}^{1+2\iota_6}$ is denoted as
$$
\a=(\a_0,\a_1,\a_{\ol1},...,\a_{\iota_6},\a_{\ol\iota_6})
\mbox{ \ with \ }\a_p\in\mbb{F}\ \ \mbox{for all}\ \ p\in\wh J.
\eqno(2.10)$$
For $\a\in\mbb{F}^{1+2\iota_6}$ and $K\subset J$, we denote by $\a_{\ssc K}$
the vector obtained from $\a$ with support $K$, i.e.,
$$
\a_{\ssc K}=(\b_0,\b_1,\b_{\ol1},...,\b_{\iota_6},\b_{\ol\iota_6})
\ \mbox{ with }\
\b_p=0\mbox{ \ if \ }p\notin K\mbox{ \ and \ }\b_p=\a_p\mbox{ \ if \ }p\in K.
\eqno(2.11)$$
Moreover, we denote
$$a_{[p]}=(0,...,0,\stl{p}{a},0,...,0)\qquad\for\;\;a\in\mbb{F}.\eqno(2.12)$$
When the context is clear, we also use the notation $\a_K$ to denote the element in
$\mbb{F}^{|K|}$ obtained from the element $\a$ in $\mbb{F}^{1+2\iota_6}$ by deleting the coordinates at $p\in \ol{0,2\iota_6}\setminus K$; for instance, $\a_{\{1,2\}}=(\a_1,\a_2).$

Take
$$\si_p=\si_{_{\ol p}}=\left\{\begin{array}{ll}
0\mbox{ or }1_{[0]}&\mbox{if \ }p=0,\vs{4pt}\\
-1_{[p]}-1_{[\ol p]}&\mbox{if \ }p\in I_{1,3},\vs{4pt}\\
-1_{[p]}&\mbox{if \ }p\in I_{4,5},\vs{4pt}\\
0&\mbox{if \ }p\in I_6.\\
\end{array}\right.
\eqno(2.13)$$
Let $\G$ be an additive subgroup of $\mbb{F}^{1+2\iota_6}$ such that
$$
\G\subset\{\a\in\mbb{F}^{1+2\iota_6}\,|\,\a_{I_6\cup\ol I_{4,6}}=0\},
\eqno(2.14)$$
$$\{\si_p\,|\,p\in\wh I\}\subset\G\mbox{ \ and \ }\mbb{F}1_{[p]}\cap\G\ne\{0\}\ \ \for\ \ p\in J_{1,3},
\eqno(2.15)$$
and
$$
\mbb{F}1_{[0]}\cap\G\ne\{0\}
\mbox{ if }\G_0\ne\{0\}.
\eqno(2.16)$$
Moreover, we define
$$\G_p=\{\a_p\mid (\a_0,\a_1,\a_{\ol1},...,\a_{\iota_6})\in\G\}.
\eqno(2.17)$$

Take
$$\JJ_0=\{0\}\mbox{ or }\mbb{N}
\mbox{ such that }\JJ_0+\G_0\ne\{0\},\eqno(2.18)$$
$$\JJ_1=\{\ii\in\mbb{N}\:^{2\iota_6}\,|\,\ii_{I_{1,2}\cup I_4\cup \ol I_1}=0\},\eqno(2.19)$$
and set
$$ \JJ=(\JJ_0,\JJ_1)\subset \mbb{N}^{1+2\ell_6},\eqno(2.20)$$
an additive subsemigroup. An element of $\JJ$  is denoted as
$$
\vec i=(i_0,\ii)=(i_0)_{[0]}+\ii,\mbox{ with }i_0\in\JJ_0,\ii\in\JJ_1.
\eqno(2.21)$$

Let $\AA$ be the semigroup algebra $\mbb{F}[\G\times\JJ]$ with a basis
$\{x^{\a,\vec i}\,|\,(\a,\vec i)\in\G\times \JJ\}$ and the multiplication $\cdot$ defined by
$$x^{\a,\vec i}\cdot x^{\b,\vec j}=x^{\a+\b,\vec i+\vec j}
\ \ \for\ \ (\a,\vec i),(\b,\vec j)\in\G\times \JJ.
\eqno(2.22)$$
Then $(\AA,\cdot)$ forms a commutative associative algebra with $1=x^{0,0}$
as the identity element. For convenience, we often denote
$$x^{\a}=x^{\a,0},\ \
t^{\vec i}=x^{0,\vec i},\ \
t_p=t^{1_{[p]}}\qquad \for\;\;\a\in\G,\;\vec i\in{\cal J},\; p\in\wh J.\eqno(2.23)$$
In particular,
$$t^{\vec i}=\prod_{p\in\wh J}t_p^{i_p}.
\eqno(2.24)$$

Define the derivations $\{\ptl_p,\ptl^*_p,\ptl_{t_p}\,|\,p\in\wh J\}$ by
$$\ptl_p=\ptl^*_p+\ptl_{t_p}\mbox{ \ and \ }
\ptl^*_p(x^{\a,\vec i})=\a_p x^{\a,\vec i},\ \
\ptl_{t_p}(x^{\a,\vec i})=i_p x^{\a,\vec i-1_{[p]}},
\eqno(2.25)$$
for $p\in\wh J,(\a,\vec i)\in\G\times\JJ,$
where we adopt the convention that if a notion is not defined but technically
appears in an expression, we always treat it as zero; for instance,
$x^{\a,-1_{[1]}}=0$ for any $\a\in\G$.
In particular,
$$\ptl^*_p=0,\ \ \ptl_{t_q}=0\ \ \for\ \ p\in
I_6\cup\ol I_{4,6},\
q\in I_{1,2}\cup I_4\cup \ol I_1,
\eqno(2.26)$$
by (2.14) and (2.19).
We call the nonzero derivations $\ptl^*_p$ {\it grading operators},
the nonzero derivations $\ptl_{t_q}$ {\it down-grading operators},
and the derivations $\ptl_r^{\ast}+\ptl_{t_r}$ {\it mixed operators}
if both $\ptl_r^{\ast}$ and $\ptl_{t_r}$ are not zero.
The types of derivation pairs in the order of the groups
$\{(\ptl_p,\ptl_{\bar{p}})\mid p\in I_i\}$ for $i\in\ol{1,6}$ are
$$
(g,g),\;(g,m),\;(m,m),\;(g,d),\;(m,d),\;(d,d),
\eqno(2.27)$$
where ``m'' stands for mixed operators, ``g'' stands for grading
operators and ``d'' stands for down-grading operators.
\par
We denote
$$
\ptl=\sum_{p\in J_{1,3}\cup I_{4,5}}
\ptl^*+\sum_{p\in I_6\cup\ol I_{4,6}}t_p\ptl_{t_p}.
\eqno(2.28)$$
Define the following Lie bracket $[\cdot,\cdot]$ on $\AA$:
$$
[u,v]=
\sum_{p\in I}x^{\si_p}
(\ptl_p(u)\ptl_{\ol p}(v)-\ptl_{\ol p}(u)\ptl_p(v))+
x^{\si_0}((2-\ptl)(u)\ptl_0(v)-\ptl_0(u)(2-\ptl)(v)),
\eqno(2.29)$$
for $u,v\in\AA$ (cf.~(2.7)). Then $(\AA,[\cdot,\cdot])$
forms a contact Lie algebra, which is in general not finitely-graded.
The algebras $(\AA,[\cdot,\cdot])$ are the normalized forms
of the contact Lie algebras constructed in [X1], which were proved to be simple
if $(\G_0.0,...,0)\subset \G$ (cf.~(2.17)).

We denote the Lie algebra $(\AA,[\cdot,\cdot])$ by
$$
{\cal K}(\vec\ell,\si,\G,\JJ),\mbox{ where }\si=\sum_{p\in\wh I_{1,5}}\si_p.
\eqno(2.30)$$

\section{Lemmas}

In this section, we shall present five lemmas as the preparation of the proof of the main Theorem.

By (2.14), (2.18), (2.19), (2.25), (2.26) and (2.28), we give a more detailed formula of  (2.29):
\begin{eqnarray*}& &
[x^{\a,\dis\vec i},x^{\b,\vec j}]
\vs{4pt}\\ &=&
\sum_{p\in I_{1,3}}(\a_p\b_{\ol p}-\a_{\ol p}\b_p)
x^{\si_p+\a+\b,\vec i+\vec j}
+\sum_{p\in I_{2,5}}(\a_pj_{\ol p}-i_{\ol p}\b_p)
x^{\si_p+\a+\b,\vec i+\vec j-1_{[\ol p]}}
\\ & & +\sum_{p\in I_3}(i_p\b_{\ol p}-j_p\a_{\ol p})
x^{\si_p+\a+\b,\vec i+\vec j-1_{[p]}}
+\sum_{p\in I_3\cup I_{5,6}}(i_pj_{\ol p}-i_{\ol p}j_p)
x^{\si_p+\a+\b,\vec i+\vec j-1_{[p]}-1_{[\ol p]}}
\\ & & +((2-\vt(\a,\vec i))\b_0-\a_0(2-\vt(\b,\vec j)))
x^{\si_0+\a+\b,\vec i+\vec j}\\ & &+((2-\vt(\a,\vec i))j_0-i_0(2-\vt(\b,\vec j)))
x^{\si_0+\a+\b,\vec i+\vec j-1_{[0]}}\hspace{4.55cm}(3.1)\end{eqnarray*}
for $(\a,\vec i),(\b,\vec j)\in\G\times\J$, where
$$
\vt(\a,\vec i)=\sum_{p\in J_{1,3}\cup I_{4,5}}\a_p+
\sum_{p\in I_6\cup\ol I_{4,6}}i_p\qquad\for\ (\a,\vec i)\in\G\times\JJ.
\eqno(3.2)$$
\par

Here is our first lemma.
\par
{\bf Lemma 3.1}. {\it For any Lie algebra ${\cal K}(\vec\ell,\si,\G,\JJ)$
of contact type with $\si_0=1_{[0]}$, there exists a Lie algebra
${\cal K}(\vec\ell,\si',\G',\JJ)$ with $\si'_0=0$ such that
${\cal K}(\vec\ell,\si,\G,\JJ)
\cong{\cal K}(\vec\ell,\si',\G',\JJ)$.
}
\par
{\it Proof.}
Define
$$\G'=\{\a+((\vt(\a,0)+n(1-\dlt_{{\cal J},\{0\}}))/2)\si_0\,|\,\a\in\G,\;n\in\mbb{Z}\},\,\si'_0=0,\,\si'_p=\si_p,\,p\in I.
\eqno(3.3)$$
Then we have another Lie algebra ${\cal K}(\vec\ell,\si',\G',\JJ)$ of contact type.
Define a bijective map $\tau:\G\times\JJ \rar\G'\times\JJ'$ such that
$$
\tau(\a,\vec i)=(\a+(\vt(\a,\vec i)/2-1)\si_0,\vec i)
\in\G'\times\JJ'\mbox{ for }(\a,\vec i)\in\G\times\JJ,
\eqno(3.4)$$
and define a linear map
$\th:{\cal K}(\vec\ell,\si,\G,\JJ)\rar{\cal K}(\vec\ell,\si',\G',\JJ)$
such that
$$
\th(x^{\a,\vec i})=x^{\tau(\a,\vec i)}\mbox{ for }(\a,\vec i)\in\G\times\JJ.
\eqno(3.5)$$
Then $\vt(\tau(\a,\vec i))=\vt(\a,\vec i)$, and by (3.1), it is straightforward
to verify that $\th([x^{\a,\vec i},x^{\b,\vec j}])=
[\th([x^{\a,\vec i}),\th(x^{\b,\vec j})]$
for $(\a,\vec i),(\b,\vec j)\in\G\times\JJ$,
i.e., $\th$ is an isomorphism.$\qquad\Box$
\psp

The above Lemma tells us that it suffices
to consider the contact Lie algebras of the form $\AA={\cal K}(\vec\ell,\si,\G,\JJ)$  with $\si_0=0$. In the rest of the paper, we always take
$$
\si_0=0.
\eqno(3.6)$$
\par
For any Lie algebra ${\cal L}$, the {\it adjoint operator} $\ad_u$ of an element $u\in {\cal L}$ is defined by
$$\ad_u(v)=[u,v]\qquad\for\;\;v\in {\cal L}.\eqno(3.7)$$
The element $u$ is called {\it ad-locally-finite} if $\ad_u$ is locally-finite (cf.~(1.23)).  Moreover, it is called {\it ad-locally nilpotent} if for any $v\in{\cal L}$, there exists an integer $n$ (depending on $v$) such that
$$(\ad_u)^n(v)=0.\eqno(3.8)$$
We denote by ${\cal L}^{\rm F}$ the set of {\it ad}-locally-finite elements in ${\cal L}$   and by ${\cal L}^{\rm N}$ the set of  {\it ad}-locally-nilpotent elements in ${\cal L}$. Furthermore, we write
$$A_0=\{1-\d_{\G_0,\{0\}},\d_{\G_0,\{0\}}t_0,x^{-\si_p},
x^{-\si_q,1_{[\ol q]}},t^{1_{[r]}+1_{[\ol r]}}
\,|\,p\in I_{1,3},q\in I_{4,5},r\in I_6\},\eqno(3.9)$$
$$A_1=\{\d_{\G_0,\{0\}},x^{\a,\vec i}\,|\,\vt(\a,\vec i)=2,\a_{\wh J_{1,3}}=
\vec i_{\wh J_{1,3}\cup\ol I_{4,5}}=0,i_pi_{\ol p}=0,\forall\,p\in I_6\},\eqno(3.10)$$
$$A_2=\{x^{\a,\vec i}\,|\,\a_{\wh J_{1,3}}=\vec i_{J_{1,3}\cup\ol I_{4,5}}=0,
(1-\d_{\G_0,\{0\}})(\vt(\a,\vec i)-2)=0,i_0\le \d_{\G_0,\{0\}}\}.\eqno(3.11)$$
For any $\vec i\in\JJ$, we define the {\it level} of $\vec i$ to be
$$
|\vec i|=\sum_{p\in\wh J}i_p.
\eqno(3.12)$$
For any $(\a,\vec i)\in\G\times\JJ$, we define the {\it support}
of $(\a,\vec i)$ to be
$$
{\rm supp}(\a,\vec i)=\{p\in\wh J\,|\,\a_p\ne0\mbox{ or }i_p\ne0\}.
\eqno(3.13)$$
\par
{\bf Lemma 3.2}. {\it (1)
$A_0\cup A_1\subset\AA^{\rm F}\subset{\it Span}(A_0\cup A_2)$. (2)
$A_1\subset\AA^{\rm N}\subset{\it Span\ssc}(A_2)$.
}
\par
{\it Proof}. Note that by (3.1), we have
$$[1,x^{\b, \vec j}]=2\b_0x^{\b,\vec j}
+2j_0x^{\b,\vec j-1_{[0]}},\eqno(3.14)$$
$$[t_0,x^{\b, \vec j}]=(2j_0-(2-\vt(\b,\vec j)))x^{\b,\vec j}\mbox{ if }\G_0=\{0\},\eqno(3.15)$$
$$[x^{-\si_p},x^{\b,\vec j}]=\left\{
\begin{array}{ll}
(\b_{\ol p}-\b_p)x^{\b,\vec j}& \for\ p\in I_1,
\vs{4pt}\\
(\b_{\ol p}-\b_p)x^{\b,\vec j}+j_{\ol p}x^{\b,\vec j-1_{[\ol p]}}
& \for\ p\in I_2,
\vs{4pt}\\
(\b_{\ol p}-\b_p)x^{\b,\vec j}-j_px^{\b,\vec j-1_{[p]}}
+j_{\ol p}x^{\b,\vec j-1_{[\ol p]}}& \for\ p\in I_3,
\end{array}
\right.\eqno(3.16)$$
$$[x^{-\si_q,1_{[\ol q]}},x^{\b,\vec j}]=\left\{
\begin{array}{ll}
(-\b_q+j_{\ol q})x^{\b,\vec j}& \for\ q\in I_4,
\vs{4pt}\\
(-\b_q+j_{\ol q})x^{\b,\vec j}
-j_qx^{\b,\vec j-1_{[q]}}
& \for\ q\in I_5,
\end{array}
\right.\eqno(3.17)$$
$$[t^{1_{[r]}+1_{[\ol r]}},x^{\b,\vec j}]=
(j_{\ol r}-j_r)x^{\b,\vec j}\ \for\ r\in I_6.\eqno(3.18)$$
Thus $A_0\subset\AA^{\rm F}$. If $\d_{\G_0,\{0\}}=1$, using (3.14) and (3.15), we see that
$\d_{\G_0,\{0\}}\in\AA^{\rm N}\subset\AA^{\rm F},t_0\in\AA^{\rm F}$. Suppose
$\d_{\G_0,\{0\}}\ne x^{\a,\vec i}\in A_1$. Then ${\rm supp}(\a,\vec i)\subset
I_{4,6}\cup\ol I_6$. Moreover,  $\ol p\notin{\rm supp}(\a,\vec i)$ if $p\in{\rm supp}(\a,\vec i)$. Let $x^{\b,\vec j}\in\AA$.
By (3.1), we see that
$[x^{\a,\vec i},x^{\b,\vec j}]$ is either zero or
a linear combination of the elements $x^{\g,\vec k}$ such that there exists at least a
$p\in (I_6\cup\ol I_{4,6})\bs{\rm supp}(\a,\vec i)$ with $k_p<j_p$.
Thus if we set
$m=1+\sum_{p\in (I_6\cup\ol I_{4,6})\bs{\rm supp}(\a,\vec i)}j_p$, then
$\ad^m_{x^{\a,\vec i}}(x^{\b,\vec j})=0$. This proves
$A_1\subset\AA^{\rm N}\subset\AA^{\rm F}$.
\par

Suppose $u\notin{\rm Span}(A_0\cup A_2)$.
Write
$$
u=\sum_{(\a,\vec i)\in M_0}c_{\a,\vec i}x^{\a,\vec i},
\eqno(3.19)$$
where
$$M_0=\{(\a,\vec i)\in\G\times\JJ\,|\,c_{\a,\vec i}\ne0\}\eqno(3.20)$$
is a finite set.
First assume that $(1-\d_{\G_0,\{0\}})(\vt(\a,\vec i)-2)=0$ for all
$(\a,\vec i)\in M_0$.
Then there exist $(\g,\vec k)\in M_0$ and
$p\in{\rm supp}(\g,\vec k)\bs J_6$ or $\ol p\in{\rm supp}(\g,\vec k)\bs J_6$ such that
$$
\begin{array}{llll}
p=0\!\!\!\!&\mbox{and}\!\!\!\!&
(\g,\vec k)\ne(0,0)\mbox{ and either }\g_0\ne0
\mbox{ or }k_0\ge \d_{\G_0,\{0\}}+1,
\!\!\!&\mbox{or}
\vs{4pt}\\
p\in I_{1,3}\!\!\!\!&\mbox{and}\!\!\!\!&
(\g,\vec k)\ne(-\si_p,0)\mbox{ and }(\g_p,\g_{\ol p},k_p,k_{\ol p})\ne0,
\!\!\!&\mbox{or}
\vs{4pt}\\
p\in I_{4,5}\!\!\!\!&\mbox{and}\!\!\!\!&
(\g,\vec k)\ne(-\si_p,1_{[\ol p]})\mbox{ and }k_{\ol p}\ne0.
\!\!\!\!&
\end{array}
\eqno(3.21)$$
Choose a total order on $\G$ compatible with group
structure of $\G$ and define the total order on $\G\times\JJ$ by
the lexicographical order,
such that the maximal element $(\g,\vec k)$ of $M_0$ satisfies
(3.21) for some $p\in{\rm supp}(\g,\vec k)\bs J_6$, and that
$\si_p>\si_q$ for all $q\ne p$. This is possible
because the set of all nonzero $\si_q$ is $\mbb{F}$-linearly independent.
Say, $p\in I_{1,3}$ and $(\g_p,\g_{\ol p})\ne0$
(the proof for other $p$ is similar).
Choose
$\b=a_{[\ol p]}\in\G\bs\{0\}$ (cf.~(2.15)) such that
$\g_pa+m(\g_{\ol p}-\g_p)\ne0$ for all $m\in\mbb{N}$.
Then
the ``highest'' term of $\ad_u^n(x^{\b})$ is
$x^{\b+n\g+n\si_p,n\vec k}$ with the coefficient
$$
\prod_{m=0}^{n-1}(\g_p(\b_{\ol p}+m\g_{\ol p}-m)-
\g_{\ol p}(m\g_p-m))=
\prod_{m=0}^{n-1}(\g_p\b_{\ol p}+m(\g_{\ol p}-\g_p))\ne0,
\eqno(3.22)$$
which implies
$$
{\rm dim}({\rm Span}\{\ad_u^n(x^{\b})\,|\,n\in\mbb{N}\})=\infty.
\eqno(3.23)$$
Thus $u\notin\AA^{\rm F}$.
Next assume that there exists $(\a,\vec i)\in M_0$ such that
$(1-\d_{\G_0,\{0\}})(\vt(\a,\vec i)-2)\ne0$. Then $\G_0\ne\{0\}$ and
$\vt(\a,\vec i)\ne2$. By making use of the last 4 terms of (3.1), we
can prove $u\notin\AA^{\rm F}$ as above.
This proves
$\AA^{\rm F}\subset{\rm Span}(A_0\cup A_2)$. Similarly,
$\AA^{\rm N}\subset{\rm Span\ssc}(A_2).\qquad\Box$

For any subset $S\subset\AA$, we define
$$
E(S)=\{u\in\AA\mid [S,u]\subset \mbb{F}u\}.
\eqno(3.24)$$
Define a map $\pi:\G\rar
\mbb{F}^{1+\iota_5}$ by
$$\pi(\a)=(\mu_0,\mu_1,...,\mu_{\iota_5})\eqno(3.25)$$
with $$
\mu_p=\left\{
\begin{array}{ll}
2\a_0&\mbox{if }p=0,\G_0\ne\{0\},
\vs{4pt}\\
\vt(\a,0)-2&\mbox{if }p=0,\G_0=\{0\},
\vs{4pt}\\
\a_{\ol p}-\a_p&\mbox{if }p\in I_{1,3},
\vs{4pt}\\
-\a_p&\mbox{if }p\in I_{4,5},
\end{array}
\right.
\eqno(3.26)$$
(cf.~(3.14)-(3.18)).
For $\mu\in\pi(\G)$, we define
$$\BB_\mu={\rm Span}\{x^\a\in\AA\,|\,
\pi(\a)=\mu\},\eqno(3.27)$$
$$\BB=\bigoplus_{\mu\in\pi(\G)}\BB_\mu
={\rm Span}\{x^\a\,|\,\a\in\G\}.\eqno(3.28)$$
Then we have
\par
{\bf Lemma 3.3}.
{\it
$$
E(\AA^{\it F})=\bigcup_{\mu\in\pi(\G)}\BB_\mu
\mbox{ (thus $\BB={\it Span\,}(E(\AA^{\it F}))$\,)}.
\eqno(3.29)$$
\par}
{\it Proof}.
We want to prove
$$
E(A_0\cup A_1)
\subset
\bigcup_{\mu\in\pi(\G)}\BB_\mu
\subset E({\rm Span}(A_0\cup A_1)).
\eqno(3.30)$$
For any $\mu\in\pi(\G)$,
elements in $\BB_\mu$ are common eigenvectors of $\ad_{A_0}$
by (3.9), (3.14)-(3.18) and (3.26)-(3.28),
and $\ad_{A_1}$ acts trivially on $\BB_\mu$.
Since elements in $A_0$ commute with each other,
elements in $\BB_\mu$ are common eigenvectors
of $\ad_{{\rm Span}(A_0\cup A_1)}$. Thus
$\cup_{\mu\in\pi(\G)}\BB_\mu
\subset E({\rm Span}(A_0\cup A_1)).$
Suppose that  $u=\sum c_{\a,\vec i}x^{\a,\vec i}\in\AA$ is a  common eigenvector of $A_0\cup A_1$. Since $\ad_{A_1}$ is locally nilpotent, $\ad_{A_1}$ acts
trivially on $u$. If $(\a,\vec i)\in M_0$ with
$i_{\ol p}\ne0$ for some $p\in I_{4,5}\cup J_6$,
then we can choose $v=x^{2\si_p}\in A_1$ if $p\in I_{4,5}$ or
$v=t^{2_{[p]}}\in A_1$ if $p\in J_6$ such that $[v,x^{\a,\vec i}]\ne0$ and
so $[v,u]\ne0$.
So we must have $\vec i_{\ol I_{4,5}\cup J_6}=0$ if $c_{\a,\vec i}\ne0$.
Similarly, since $u$ is a common eigenvector of $\ad_{A_1}$, we must
have by (3.14)-(3.18) that $\vec i_{\wh J_{1,3}\cup I_5}=0$ (thus $\vec i=0$)
and $\pi(\a)=\mu$ for some $\mu$ if $c_{\a,\vec i}\ne0$.
This shows that $u\in\BB_\mu$. Hence (3.30) is proved.
Using Lemma 3.2(1), we have
$$
E(A_0\cup A_1)\supset E(\AA^{\rm F})\supset
E({\rm Span}(A_0\cup A_1)).\eqno(3.31)$$
This and (3.30) show that all these sets are equal, i.e., we have (3.29).$\qquad\Box$

{\bf Lemma 3.4}.
{\it Denote
$$
\BB_{F}=
\{x^{-\si_p},x^\a\,|\,p\in\wh I_{1,3},(1-\d_{\G_0,\{0\}})(\vt(\a,0)-2)=0,
\a=\a_{I_{4,5}}\},
\eqno(3.32)$$
$$
\BB_{N}=
\{x^\a\,|\,(1-\d_{\G_0,\{0\}})(\vt(\a,0)-2)=0,\a=\a_{I_{4,5}}\}.
\eqno(3.33)$$
Then $\BB^{\rm F}=\mbox{\it Span}(\BB_{\rm F}),\,\BB^{\rm N}=\mbox{\it Span}(\BB_{\rm N})$.
\par}
{\it Proof.}
We shall prove that $\BB^{\rm F}={\rm Span}(\BB_{\rm F})$ as the proof
that $\BB^{\rm N}={\rm Span}(\BB_{\rm N})$ is similar.
It is straightforward to verify that elements in $\BB_{\rm F}$ are
{\it ad}-locally finite on $\BB$ and that
$\BB_{\rm F}$ is commutative.
Thus ${\rm Span}(\BB_{\rm F})\subset\BB^{\rm F}$.
Conversely, suppose $u\in\BB\bs{\rm Span}(\BB_{\rm F})$. Then we can write $u$
as in (3.19),
where now
$$M_0=\{(\a,\vec i)\in\G\times\JJ_{\wh J_{1,5}}\,|\,
\vec i=0,c_{\a,\vec i}\ne0\}\eqno(3.34)$$
is finite.
Thus we still have (3.21), and the same arguments after (3.21) show
that $u$ is not {\it ad}-locally finite on $\BB.\qquad\Box$

We also have that the center of $\BB$ is
$$
C(\BB)=
{\rm Span}\{x^\a\,|\,
(1-\d_{\G_0,\{0\}})(\vt(\a,0)-2)=0,\a=\a_{I_{4,5}}\in\G\}
\eqno(3.35)$$
by (3.1) and (3.28).

\par
{\bf Lemma 3.5}. {\it
(1)
Assume that $\G_0\ne\{0\}$.
For $\mu\in\pi(\G)$, regarding
$\BB_\mu$ as a $\BB_0$-module, we have (i)
if $\mu_{\wh I_{1,3}}=0$, the action of $\BB_0$ on $\BB_\mu$
is trivial and (ii) if $\mu_{\wh I_{1,3}}\ne0$,
$\BB_\mu$ is a cyclic $\BB_0$-module, and
the nonzero scalar multiples of $x^{\a}$ for all $\a\in\G$ with
$\pi(\a)=\mu$ are the only generators.
\par
(2) Assume that $\G_0=\{0\}$ and $\iota_3\ne0$.
We have (i)
if $\mu_{I_{1,3}}=0$, the action of $\BB_0$ on $\BB_\mu$
is trivial and (ii) if $\mu_{I_{1,3}}\ne0$,
$\BB_\mu$ is a cyclic $\BB_0$-module, and
the nonzero scalar multiple of $x^{\a}$ for all $\a\in\G$ with
$\pi(\a)=\mu$  are the only generators.
\par
(3) Assume that $\G_0=\{0\},\iota_3=0$ and $\ell_4+\ell_5\ne0$.
Then $(\cup_{\a\in\G}\mbb{F} x^\a)\bs\{0\}$ is the set of the common eigenvectors
of $\AA^{\rm F}$ in $\BB$.
}
\par
{\it Proof}. Note that in $\BB$, we have
$$
[x^{\a},x^{\b}]
=\sum_{p\in I_{1,3}}(\a_p\b_{\ol p}-\a_{\ol p}\b_p)
x^{\si_p+\a+\b}
+((2-\vt(\a,0))\b_0-\a_0(2-\vt(\b,0)))
x^{\a+\b}\eqno(3.36)$$
for $\a,\b\in\G$ by (3.1).
\par
(1)
Assume that $\G_0\ne\{0\}$.
One can easily verify that
$$
\pi(\si_p+\a+\b)=\pi(\a)
+\pi(\b)\ \ \ \for\ \ \ p\in\wh I_{1,3}
\eqno(3.37)$$
by (3.25)-(3.26).
{}From (3.36)
we see that $x^{\a}$ commutes with
$x^{\b}$ if $\pi(\a)=0$ and $(\pi(\b))_{\wh I_{1,3}}=0$
(cf.~(3.25)-(3.26)). Thus
if $\mu_{\wh I_{1,3}}=0$, then the adjoint action of $\BB_0$ on $\BB_\mu$
is trivial.
Assume that $u=\sum_{\b\in M_0}c_{\b}x^{\b}\in\BB_\mu$ with
$\mu_{\wh I_{1,3}}\ne0$, where $M_0=\{\b\in\G\,|\,\pi(\b)=\mu,
c_{\b}\ne0\}$.
By (3.36), one can verify that
$$
[x^{\a},u]=(\sum_{p\in I_{1,3}}\a_p\mu_p x^{\si_p+\a}
+(1-\vt(\a,0)/2)\mu_0x^{\a})u\mbox{ if }\pi(\a)=0,
\eqno(3.38)$$
where by (3.25)-(3.26), the coefficient $(1-\vt(\a,0)/2)\mu_0$ comes from
$$
(2-\vt(\a,0))\b_0-\a_0(2-\vt(\b,0))=
(2-\vt(\a,0))\b_0=(1-\vt(\a,0)/2)\mu_0.
\eqno(3.39)$$
Thus
$$
U={\rm Span}\{x^{\si_p+\a}u=
\sum_{\b\in M_0}c_{\b}x^{\si_p+\a+\b}\,
|\,\a\in\kn_{\pi},p\in\wh I_{1,3}\}
\eqno(3.40)$$
is a $\BB_0$-submodule of $\BB_\mu$.

Let $\la u\ra$ denote
the cyclic submodule of $\BB_\mu$ generated by $u$. Then
$\la u\ra\subset U$. If the size $|M_0|$ of $M_0$ is $\ge2$, then
$U$ is a proper submodule of $\BB_\mu$, and so $u$ is not a generator.
Suppose $M_0$ is a singleton $\{\b\}$ with $\pi(\b)=\mu$.
Say $\mu_0\ne0$ (the proof is similar if $\mu_p\ne0$ for some $p\in I_{1,3}$).
If $\iota_3=0$, then $\BB_\mu$ is 1-dimensional and
we see that the proof is trivial. Assume that $\iota_3>0$.
For $p\in I_{1,3},k,m\in\mbb{Z}$, take $\g_k=k\si_p$. Then
by (3.25)-(3.26), $\g_k\in\kn_{\pi}$.
By (3.36), we have
$$[x^{\g_k},x^{\b}]=-k\mu_px^{\b+(k+1)\si_p}+
(1+k)\mu_0x^{\b+k\si_p}\in\la u\ra,\eqno(3.41)$$
\begin{eqnarray*}& &[x^{\g_{m-k}},[x^{\g_k},x^{\b}]]
\\ &=&-k\mu_p(-(m-k)\mu_p x^{\b+(m+2)\si_p}+(1+m-k)\mu_0x^{\b+(m+1)\si_p})
\\ & &+(1+k)\mu_0(-(m-k)\mu_p x^{\b+(m+1)\si_p}
+(1+m-k)\mu_0x^{\b+m\si_p})\in\la u\ra.\hspace{2.25cm}(3.42)\end{eqnarray*}
{}From (3.41) and (3.42), one deduces that
$$
x^{\b+m\si_p}\in\la u\ra\qquad\for\ \ m\in\mbb{Z}.
\eqno(3.43)$$
For $\a\in\kn_{\pi}$, we see that
$\a-(k+1)\si_p\in\kn_{\pi}$, and
by (3.36),
\begin{eqnarray*}& &[x^{\a-k\si_p},x^{\b+k\si_p}]\\
&=&{\dis\sum_{q\in I_{1,3}}}(\a_q+k\d_{p,q})\mu_qx^{\si_q+\a+\b}
+(1-k-\vt(\a,0)/2)\mu_0x^{\a+\b}\in\la u\ra.\hspace{2.55cm}(3.44)\end{eqnarray*}
Comparing the coefficient of $k$, we obtain
$$
u_p=\mu_px^{\si_p+\a+\b}-\mu_0x^{\a+\b}\in\la u\ra
\qquad\for\;\; p\in I_{1,3}.
\eqno(3.45)$$
Note that since $\pi(\a)=0$,  (3.2) and (3.25)-(3.26) imply
$$\vt(\a,0)=\sum_{p\in J_{1,3}\cup I_{4,5}}\a_p=2\sum_{p\in I_{1,3}}\a_p,\eqno(3.46)$$
By (3.45),  the right-hand side of (3.44) becomes
$$\sum_{p\in I_{1,3}}u_p+\mu_0x^{\a+\b}.\eqno(3.47)$$
This shows that $x^{\a+\b}\in\la u\ra$ for all $\a
\in\kn_{\pi}$. Since $\BB_\mu$ is spanned by such elements,  $u$ is a generator of
$\BB_\mu$.

The proof of (2) is similar. (3) is obvious by (3.29).$\qquad\Box$

\section{Main Theorem}

In this section, we shall present the main theorem on isomorphism classes of the contact Lie algebras.

We assume that $\mbb{F}$ is algebraically closed in this section.
Denote by $M_{m\times n}(\mbb{F})$ the space of $m\times n$ matrices with entries
in $\mbb{F}$ and by $GL_m(\mbb{F})$ the group of $m\times m$ invertible matrices with
entries in $\mbb{F}$. Define
$$
G_p=\left\{\left(\begin{array}{cc}
1-a&1-a-b\\ a&a+b\end{array}\right)\mid a,b\in\mbb{F},\;b\ne0\right\}
\qquad \for\;\;p\in I_1\cup I_3,
\eqno(4.1)$$
$$
G_q=\left\{\left(\begin{array}{cc}
1&1-b\\ 0&b\end{array}\right)\mid 0\neq b\in\mbb{F}\right\}
\qquad \for\;\;q\in I_2,
\eqno(4.2)$$
and denote
$$
G=\{{\rm diag}(b_0,g_1,...,g_{\iota_3})\in GL_{1+2\iota_3}(\mbb{F})
\mid 0\ne b_0\in\mbb{F},\,g_i\in G_i\}.
\eqno(4.3)$$
Define $H$ to be the set of $(1+2\iota_3)\times(\ell_4+\ell_5)$ matrices of
the form
$$
h=\left(
\begin{array}{l}h_0\\ h_1\\ 0\\ \end{array}\right),
\mbox{ where }h_0\in M_{1\times(\ell_4+\ell_5)}(\mbb{F}),\,
h_1\in M_{2\ell_1\times(\ell_4+\ell_5)}(\mbb{F})\mbox{ and}
\eqno(4.4)$$
$$
h_0=0\mbox{ if }J_0=\mbb{N}\mbox{ and }h_1=(h_{p,q})\mbox{ with }h_{2p,q}=-h_{2p-1,q}
\mbox{ for }p\in I_1.
\eqno(4.5)$$
Define $F$ to be the group of invertible matrices of the form
$$
f=\left(\begin{array}{ccc} A&C\\ 0&B
\end{array}\right),\mbox{ where }A\in GL_{\ell_4}(\mbb{F}),\,B\in GL_{\ell_5}(\mbb{F}),\,
C\in M_{\ell_4\times\ell_5}(\mbb{F}).
\eqno(4.6)$$
Define a subgroup of $GL_{1+2\iota_6}$,
$$
\begin{array}{l}
{\cal G}''=GL_{1+2\iota_6}\mbox{ if }\G_0=\{0\},\mbox{ or otherwise}
\vs{4pt}\\ \dis
{\cal G}''=
\{(a_{p,q})\in GL_{1+2\iota_6}(\mbb{F})\,|\,
\sum_{q=2}^{2\iota_6+1}a_{1,q}=0,\sum_{q=1}^{2\iota_6+1}a_{p,q}=1,
a_{p,1}=0
\vs{4pt}\\ \hskip 5cm
\mbox{ for }p\in\ol{2,2\iota_6+1}\}.
\end{array}
\eqno(4.7)$$
Denote by ${\bf 1}_n$ the $n\times n$ identity matrix.
Now we define the group
$$
{\cal G}'=\left\{
\left(\begin{array}{ccc}
g&h&0\\ 0&f&0\\ 0&0&{\bf1}_{\ell_4+\ell_5+2\ell_6}\end{array}\right)
\in {\cal G}''\mid g\in G,\,h\in H,\,f\in F\right\}.
\eqno(4.8)$$

Let ${\cal S}_{\wh I}$ be the permutation group on the index set $\wh I$
(cf.~(2.5) and (2.6)). Define the subgroup
$$
\begin{array}{ll}
{\cal S}=\{\nu\in {\cal S}_{\wh I}\;|\!\!\!\!&
\nu(\wh I_1)=\wh I_1,\,\nu(I_2)=I_2,\,\nu(I_3)=I_3,\,
\nu|_{I_{4,6}}=\mbox{Id}_{I_{4,6}}\mbox{ and}
\vs{4pt}\\ &
\nu(0)=0\mbox{ if }J_0=\mbb{N}\}.
\\
\end{array}
\eqno(4.9)$$
For $\nu\in{\cal S},\a\in\mbb{F}^{1+2\iota_6}$, we define
$$
\nu(\a)=\left\{
\begin{array}{l}
(\a_{0},\a_{\nu(1)},\a_{\ol{\nu(1)}},...,\a_{\nu(\iota_6)},
\a_{\ol{\nu(\iota_6)}})\mbox{ if }\nu(0)=0,
\vs{4pt}\\
\nu'(\a)\mbox{ otherwise,}
\\
\end{array}
\right.\eqno(4.10)$$
where
$$\nu'\in{\cal S}\mbox{ such that }\nu'(0)=0,\,\nu'(\nu^{-1}(0))=\nu(0),\,
\nu'|_{I\bs\{\nu^{-1}(0)\}}=\nu|_{I\bs\{\nu^{-1}(0)\}}.
\eqno(4.11)$$
Moreover, for $g'\in{\cal G}'$ and $\nu\in{\cal S}$, we define the
group automorphism $g'_{\nu}$ of $\mbb{F}^{1+2\iota_6}$ by
$$
\nu(g'_\nu(\a))=
(\a_{_{\wh J_{1,3}}},\a_{_{I_{4,5}}},\a_{_{\ol I_{4,5}\cup J_6}})g'
\mbox{ if }\nu(0)=0,\mbox{ and }
\eqno(4.12)$$
\begin{eqnarray*}\nu(g'_\nu(\a))&=&({1\over2}(\a_{\ol p}-\a_p),\a_1,\a_{\ol1},...,
\a_{p-1},\a_{\ol{p-1}},-{1\over2}\vt(\a,0)-\a_0,-{1\over2}\vt(\a,0)+\a_0,
\\ &&\a_{p+1},\a_{\ol{p+1}},...,\a_{\iota_3},\a_{\ol\iota_3},
\a_{_{I_{4,5}}},\a_{_{\ol I_{4,5}\cup J_6}})g'
\mbox{ if }p=\nu^{-1}(0)\in I_1,\hspace{2.3cm}(4.13)\end{eqnarray*}
where $\a\in \mbb{F}^{1+2\iota_6}$ and the multiplication
in the above is the vector-matrix multiplication. Define
$$
{\cal G}=\{g'_{\nu}\mid g'\in {\cal G}',\;\nu\in{\cal S}\}.
\eqno(4.14)$$ Then $\cal G$ is a subgroup of additive automorphism
of $\mbb{F}^{1+2\iota_6}$.
\par
Let ${\cal K}(\vec\ell',\si',\G',\JJ')$ be another Lie algebra defined in
Section 2. We shall add a prime on all the constructional ingredients
related to ${\cal K}(\vec\ell',\si',\G',\JJ')$ except that we use the same symbol $x$; for instance, ${\cal A}',\;\ell'_i,\;\iota'_i,$ etc.
\par
{\bf Theorem 4.1}. {\it The two Lie algebras ${\cal K}(\vec\ell,\si,\G,\JJ)$ and
${\cal K}(\vec\ell',\si',\G',\JJ')$ are isomorphic if and only if
$(\vec\ell,\JJ)=(\vec\ell',\JJ')$ and there exists $\tau\in{\cal G}$
such that $\G'=\tau(\G)$.}
\par
{\it Proof}.
``$\Longleftarrow$'' First we prove the sufficiency. Assume that
$(\vec\ell,\JJ)=(\vec\ell',\JJ')$ and there exist a group isomorphism
$\tau=g'_\nu\in{\cal G}$ such that $\G'=\tau(\G)$. We shall define an isomorphism $\th: {\cal K}(\vec\ell,\si,\G,\JJ)\rightarrow{\cal K}(\vec\ell',\si',\G',\JJ')$. To do this, we shall first find the images of some elements (see (4.30)-(4.32)). Denote
$$P={\rm diag}(2,(^{-1}_{{\ssc\ }1}),...,(^{-1}_{{\ssc\ }1}))\in
M_{(1+2\iota_3)\times(1+\iota_3)},\eqno(4.15)$$
$$P'={\rm diag}(P,-{\bf1}_{\ell_4+\ell_5})\in
M_{(1+\iota_3+\iota_5)\times(1+\iota_5)}.\eqno(4.16)$$
By (4.8),
$$
g'=
\left(\begin{array}{ccc}
g&h&0\\ 0&f&0\\ 0&0&{\bf1}_{\ell_4+\ell_5+2\ell_6}\end{array}\right),
\mbox{ where }g={\rm diag}(b_0,g_1,...,g_{\iota_3})\mbox{ with}
\eqno(4.17)$$
$$g_p=\left(\begin{array}{cc}1-a_p&1-a_p-b_p\\ a_p&a_p+b_p\end{array}\right),\;\;
g_q=\left(\begin{array}{cc}1&1-b_p\\ 0&b_p\end{array}\right)
\eqno(4.18)$$
for $p\in I_1\cup I_3,\;q\in I_2$, and
$h$ is of the form (4.4) and $f$ is of the form (4.6).
We shall take $b_0=1$ if $\G_0=\{0\}$ without loss of generality.
For convenience, we write $g'$ as $g'=(g'_{p,q})$ where $p,q\in\wh J$.
Denote
$$c_p=\sum_{q\in\wh J}g'_{p,q}\;\;\;\for\;\;p\in\wh J.
\eqno(4.19)$$

Note that if $\G_0\ne\{0\}$, we have $c_p=1$ for $p\in J$.
In general, we have
$$c_p-1=-(c_{\ol p}-1)\;\;\;\for\;\;p\in I_{1,3}.
\eqno(4.20)$$
Denote
$$D=(0,-c_1+1,-c_2+1,...,-c_{\iota_5}+1)^T\in M_{(1+\iota_5)\times 1},
\eqno(4.21)$$
where the up-index ``T'' means the transpose of the matrix.
By (4.4) and (4.5), we can write $h$ as the form
$$h=Ph'\qquad\mbox{for some}\;\;h'\in M_{(1+\iota_3)\times(\ell_4+\ell_5)}
\eqno(4.22)$$
(cf.~(4.15)). Also we have
$$\left(\begin{array}{cc}g&h\\ 0&f\\ \end{array}\right)P'=
P'\wt g'\mbox{ for some }\wt g'\in GL_{1+\iota_5}
\eqno(4.23)$$
(cf.~(4.16)). Denote
$$\wt g={\rm diag}(\wt g_1,...,\wt g_{\iota_3}),
\mbox{ where }\wt g_p=b_pg^{-1}_p.
\eqno(4.24)$$
For convenience, we denote
$$\ol t=(t_0,-t_{\ol1},t_1,-t_{\ol2},t_2,...,-t_{\ol\iota_6},t_{\iota_6}).
\eqno(4.25)$$
For a subset $K$ of $\wh J$, we denote by $\ol t_{_K}$ the vector obtained
from $\ol t$ by deleting $-t_{\ol p}, t_q$ for $\ol p,q\in J\setminus K$.
For instance,
$$\ol t_{\{\ol 1,\ol 2,2\}}=(-t_{\ol 1},-t_{\ol 2},t_2).
\eqno(4.26)$$
Denote the vectors
$$x^{-\si}=(1,x^{-\si_1},...,x^{-\si_{\iota_3}},
x^{-\si_{\iota_3+1},1_{[\ol{\iota_3+1}]}},...,x^{-\si_{\iota_5},
1_{[\ol\iota_5]}})\eqno(4.27)$$
(cf.~(3.6)) and
$$\vsi=(t_0,t_1,t_{\ol1},...,t_{\iota_3},t_{\ol\iota_3},
t_{\iota_3+1},x^{-\si_{\iota_3+1},1_{[\ol{\iota_3+1}]}},...,
t_{\iota_5},x^{-\si_{\iota_5},1_{[\ol\iota_5]}},
t_{\iota_5+1},t_{\ol{\iota_5+1}},...,t_{\iota_6},t_{\ol\iota_6}).
\eqno(4.28)$$
\par
{\it Case a}: First assume that $\nu(0)=0$.
Note that (4.12) guarantees
$$\tau(\si_p)=\si_{\nu(p)}\;\;\for\;\;p\in\wh I_{1,3}.
\eqno(4.29)$$
We shall find the image of $\zeta$. To do this, in ${\cal A}'$, we define
$$s_0=b_0t'_0+\nu(x^{-\si})\wt D,\;
\ol s_{J_{1,3}}=\nu(\ol t'_{J_{1,3}})\wt g+E,
\eqno(4.30)$$
$$s_{\ol I_{4,5}}=(x^{-\si'})_{I_{4,5}}f^{-1}
+\nu((x^{-\si'})_{\wh I_{1,3}})\wt h,
\eqno(4.31)$$
$$s_{I_4\cup I_6}=t'_{I_4\cup I_6},\;\;
s_{I_5}=t'_{I_5}B^T,\;\;
s_{\ol I_6}=b_0t'_{\ol I_6},
\eqno(4.32)$$
(cf.~(4.10)), where
$$\wt D=(\wt g')^{-1}D,\;
\wt h=-{\rm diag}(b_0,b_1,...,b_{\iota_3})^{-1}h'f^{-1},
\eqno(4.33)$$
(cf.~(4.22) and (4.23)) and where $E$ will be determined later in (4.56). We shall define
$s$ to be the image of $\zeta$. Expressions (4.30)-(4.32) are motivated from (4.43).
\par
Let
$$
\Dlt_0=\sum_{p=1}^{\iota_5}\mbb{Z}\si_p
\eqno(4.34)$$
be the subgroup of $\G$ generated by
$\{\si_i\mid i\in\ol{1,\iota_5}\}$ (cf.~(2.15)). Define
$\chi:\Dlt_0\rar\mbb{F}^{\times}$ to be the homomorphism from additive group
to the multiplication group of nonzero elements of $\mbb{F}$ determined by
$$\chi(\si_p)=b_0^{-1}b_p,\;\chi(\si_q)=b_0^{-1}\ \ \for\ \ p\in I_{1,3},\;q\in I_{4,5}.
\eqno(4.35)$$
We want to prove that $\chi$ can be extended to a homomorphism
$\chi:\G\rar\mbb{F}^{\times}$. Suppose that $\Dlt$ is a maximal subgroup of
$\G$ containing $\Dlt_0$ such that $\chi$ can be extended to a homomorphism
$\chi:\Dlt\rar\mbb{F}^{\times}$. Assume $\Dlt\ne\G$. We take an element
$\a\in\G\setminus\Dlt$. Set
$$\Dlt'=\mbb{Z}\a+\Dlt.
\eqno(4.36)$$
If $\mbb{Z}\a\cap\Dlt=\{0\}$, then we extend $\chi$ by
$$
\chi(m\a+\b)=\chi(\b)\qquad\for\;\;m\in\mbb{Z},\;\b\in\Dlt.
\eqno(4.37)$$
If  $\mbb{Z}\a\cap\Dlt=\mbb{Z} n\a$, we take an $n$th root $a$ of
$\chi(n\a)$ (recall that $\mbb{F}$ is algebraically closed) and extend
$\chi$ by
$$
\chi(m\a+\b)=a^m\chi(\b)\qquad\for\;\;m\in\mbb{Z},\;\b\in\Dlt.
\eqno(4.38)$$
It is straightforward to verify that $\chi: \Dlt'\rar \mbb{F}^{\times}$ is a
group homomorphism. This leads to a contradiction with the maximality of $\Dlt$.
So $\chi$ can be extended to a homomorphism $\chi: \G\rar \mbb{F}^{\times}$. Take
any such extension.
\par
Noting that $x^{\a,\vec i}=x^{\a'}\vsi^{\vec i}$, where
$\a'_p=\a_p,\a'_q=\a_q-i_{\ol q}$ for $p\in\wh J\bs I_{4,5},q\in I_{4,5}$,
we define the linear space  isomorphism $\th: {\cal A}\rar {\cal A}'$ by
$$
\th(x^\al\vsi^{\vec i})=b^{-1}_0\chi(\al)
x^{\tau(\al)}s^{\vec i}
\qquad \for\;\;(\al,\vec i)\in \G\times{\cal J}.
\eqno(4.39)$$
\par
{\bf Claim 1}. $\th$ is a Lie algebra isomorphism, i.e.,
$$
\th([x^\a\vsi^{\vec i},x^\b\vsi^{\vec j}])=[\th(x^\a\vsi^{\vec i}),
\th(x^\b\vsi^{\vec j})]
\qquad\for\;\;(\a,\vec i),(\b,\vec j)\in\G\times\JJ.
\eqno(4.40)$$
\par
{\it Case a.1}: First suppose $b_0=1$.
{}From (2.29), we have
$$[u,v\cdot w]=[u,v]\cdot w+v\cdot[u,w]+2\ptl_0(u)\cdot v\cdot w
\;\;\for\;\;u,v,w\in\AA.\eqno(4.41)$$
By (4.39), we have
$$\th(u\cdot v)=\th(u)\cdot\th(v),\;\;\th(\ptl_0(u))=\ptl_0(\th(u))
\;\;\for\;\;u,v\in\AA.
\eqno(4.42)$$
Thus it suffices to verify
$$\th([x^\a,x^\b])=[\th(x^\a),\th(x^\b)],\;
\th([\vsi_p,x^\a])=[\th(\vsi_p),\th(x^\a)],\;
\th([\vsi_p,\vsi_q])=[\th(\vsi_p),\th(\vsi_q)],
\eqno(4.43)$$
for all $\a,\b\in\G,\,p,q\in\wh I_3\cup I_{5,6}\cup \ol I_{2,6}$
and $p,q\ne0$ if $J_0=\{0\}$.
\par
By (2.10) and (3.36), we have
$$[x^{\al},x^{\be}]=\sum_{p\in I_{1,3}}
\left|\begin{array}{c}\al_{\{p,\ol p\}}\\\be_{\{p,\ol p\}}
\end{array}\right|x^{\si_p+\al+\be}+
((2-\vt(\a,0))\b_0-\a_0(2-\vt(\b,0)))x^{\al+\be},
\eqno(4.44)$$
where
$$
\left|\begin{array}{c}\al_{\{p,\ol p\}}\\\be_{\{p,\ol p\}}\end{array}\right|
=\left|\begin{array}{cc}\al_p&\al_{_{\ol p}}\\\be_p&\be_{_{\ol p}}\end{array}
\right|\eqno(4.45)$$
is a $2\times2$ determinant. Moreover, by (3.36), (4.12) and (4.39), we get
\begin{eqnarray*}[\th(x^{\al}),\th(x^{\be})]&\!\!\!\!=\!\!\!\!&
\chi(\al)\chi(\be)(\sum_{p\in I_{1,3}}
\left|\begin{array}{c}(\tau(\al))_{\{p,\bar{p}\}}\\
 (\tau(\be))_{\{p,\ol p\}}\end{array}\right|
{x'}^{\si'_{\nu(p)}+\tau(\al)+\tau(\be)}\vs{4pt}\\ &&+
((2-\vt(\tau(\a),0))(\tau(\b))_0-(\tau(\a))_0(2-\vt(\tau(\b),0)))
{x'}^{\tau(\al)+\tau(\be)})\vs{4pt}\\ &\!\!\!\!=\!\!\!\!&\th(
\sum_{p\in I_{1,3}}
\left|\begin{array}{c}\al_{\{p,\ol p\}}\\ \be_{\{p,\ol p\}}
\end{array}\right|x^{\si_p+\al+\be}
+((2-\vt(\a,0))\b_0-\a_0(2-\vt(\b,0)))x^{\a+\b})\vs{4pt}\\
&\!\!\!\!=\!\!\!\!&\th([x^\a,x^\b])
\hspace{10cm}
(4.46)\end{eqnarray*}
because (4.18) implies that the determinant of $g_p$ is $|g_p|=b_p=\chi(\si_p)$
(cf.~(4.35)) and
\begin{eqnarray*}\hspace{2cm}
\chi(\al)\chi(\be)\left|\begin{array}{c}(\tau(\al))_{\{p,\bar{p}\}}\\
(\tau(\be))_{\{p,\ol p\}}\end{array}\right|&=&
\chi(\al+\be)\left|\begin{array}{c}\al_{\{p,\ol p\}}\\
\be_{\{p,\ol p\}}\end{array}\right|\cdot|h_p|\vs{4pt}\\&=&
\chi(\si_p+\al+\be)\left|\begin{array}{c}\al_{\{p,\ol p\}}\\
\be_{\{p,\ol p\}}\end{array}\right|,
\hspace{3.1cm}
(4.47)\end{eqnarray*}
and (4.7), (4.12) imply
$$
\vt(\tau(\a),0)=\vt(\a,0)\mbox{ if }\G_0\ne\{0\},\mbox{ and }
(\tau(\a))_0=b_0\a_0=\a_0.
\eqno(4.48)$$
This proves the first equation of (4.43).

\par
Recall the notations (2.11) and (4.26), and ${\cal J}\subset\mbb{F}^{1+2\iota_6}$ (cf.~(2.18)-(2.20)).
As $1\times2$ matrices,
$$
[\ol t_{\{p,\ol p\}},x^{\al}]=([-t_{_{\ol p}},x^{\al}],[t_p,x^{\al}])=
\al_{\{p,\ol p\}}x^{\si_p+\al}\qquad\for\ \ p\in I.
\eqno(4.49)$$
Now we verify the second equation in (4.43). If $\G_0\ne\{0\}$,
we have
\begin{eqnarray*}\hspace{2cm}[\th(t_0),\th(x^\a)]&=&
\chi(\a)[s_0,x^{\tau(\a)}]=\chi(\a)[t'_0,x^{\tau(\a)}]\\
&=&\chi(\a)(2(\tau(\a))_0x^{\a,1_{[0]}}-(2-\vt(\tau(\a),0))x^{\tau(\a)})
\\ &=&\th([t_0,x^\a])\hspace{8cm}(4.50)\end{eqnarray*}
by (3.1) and (4.48).
If $\G_0=\{0\}$, we have
\begin{eqnarray*}\qquad [\th(t_0),\th(x^\a)]
&=&\chi(\a)[t'_0+\nu(x^{-\si})\wt D,x^{\tau(\a)}]
\\ &=&
\chi(\a)(-2+\vt(\tau(\a),0)+\nu(\tau(\a))_{\wh J_{1,3}\cup I_{4,5}}P'\wt D)
x^{\tau(\a)}\\ &=&\chi(\a)(-2+\vt(\a,0))x^{\tau(\a)} = \th([t_0,x^\a]),\hspace{4.3cm}(4.51)\end{eqnarray*}
(cf.~(4.16), (4.21), (4.23), (4.30) and (4.33)) by noting that
\begin{eqnarray*}\qquad\nu(\tau(\a))_{\wh J_{1,3}\cup I_{4,5}}P'\wt D
&=&\a_{\wh J_{1,3}\cup I_{4,5}}\left(
\begin{array}{ll}g&h\\ 0&f\\ \end{array}\right)P'\wt D=\a_{\wh J_{1,3}\cup I_{4,5}}P'\wt g'\wt D\\ &=&\a_{\wh J_{1,3}\cup I_{4,5}}P'D
\ =\ \vt(\tau(\a),0)-\vt(\a,0),\hspace{3cm}(4.52)\end{eqnarray*}
where the last equality follows from the definition of $\vt(\a,0)$ and
by (4.12). If $p\in I_2$, by (4.30), we have $s_{\ol p}=b_pt'_{\ol{\nu(p)}}$
if $\G_0\ne\{0\}$ or $s_{\ol p}\equiv b_pt'_{\ol{\nu(p)}}$ (mod\,$\mbb{F}$)
but $[\mbb{F},x^\a]=0$ if $\G_0=\{0\}$. Thus
\begin{eqnarray*}\qquad[\th(t_{\ol p}),\th(x^\a)]
&=&\chi(\a)[s_{\ol p},x^{\tau(\a)}]\\ &=&b_p\chi(\a)[t'_{\ol{\nu(p)}},x^{\tau(\a)}]\\ &=&
-\chi(\si_p+\a)(\tau(\a))_{\nu(p)}x^{\si_{\nu(p)}+\tau(\a)}\\
&=&
-\chi(\si_p+\a)\a_px^{\si_{\nu(p)}+\tau(\a)}
\ =\ \th([t_{\ol p},x^\a]).\hspace{4cm}(4.53)\end{eqnarray*}
Similarly, for $p\in I_3$, we have
\begin{eqnarray*}\hspace{2cm}
[\th(\ol t_{\{p,\ol p\}}),\th(x^{\al})]
&=&\chi(\al)[\ol s_{\{p,\ol p\}},{x'}^{\tau(\al)}]\vs{4pt}\\
&=&\chi(\al)[\ol t'_{\{\nu(p),\ol {\nu(p)}\}},{x'}^{\tau(\al)}]
b_pg_p^{-1}\vs{4pt}\\
&=&\chi(\al)b_p(\tau(\al))_{\{\nu(p),\ol{\nu(p)}\}}
g^{-1}_p{x'}^{\si'_{\nu(p)}+\tau(\al)}\vs{4pt}\\
&=&\chi(\al)b_p\al_{\{p,\ol p\}} {x'}^{\si'_{\nu(p)}+\tau(\al)}\vs{4pt}\\
&=& \chi(\al+\si_p)\al_{\{p,\ol p\}} {x'}^{\si'_{\nu(p)}+\tau(\al)}\vs{4pt}\\
&=&\th(\al_{\{p,\ol p\}} x^{\al+\si_p})\vs{4pt}\\
&=&\th([\ol t_{\{p,\ol{p}\}},x^{\al}]).
\hspace{6.7cm}(4.54)\end{eqnarray*}
Furthermore, $[t_{I_5\cup J_6},x^\a]=0$ for $\al\in\G$. We have
$gP=P{\rm diag}(b_0,...,b_{\iota_3})$ by (4.15), (4.17) and (4.18), and
\begin{eqnarray*}\hspace{2.5cm}
[\th(\vsi_{\ol I_{4,5}}),\th(x^\a)]&=&
\chi(\a)[s_{\ol I_{4,5}},x^{\tau(\a)}]
\\ &=&\chi(\a)[(x^{-\si'})_{I_{4,5}}f^{-1}
+\nu((x^{-\si'})_{\wh I_{1,3}})\wt h,x^{\tau(\a)}]\\ &=&
\chi(\a)((\tau(\a))_{I_{4,5}}f^{-1}+\a_{\wh J_{1,3}}gP\wt h)x^{\tau(\a)}
\\ &=&\chi(\a)\a_{I_{4,5}}x^{\tau(\a)}\\ &=&\th(\a_{I_{4,5}}x^\a)\\
&=&\th([\vsi_{\ol I_{4,5}},x^\a])\hspace{6.6cm}(4.55)\end{eqnarray*}
by (4.12), (4.22), (4.27), (4.31) and (4.33).
This proves the second equation of (4.43).
\par
Now we take $E=0$ if $\G_0\ne\{0\}$ and in general take
$$
E=2^{-1}\wt D^T[\nu(x^{-\si})^T,\nu(\ol t'_{J_{1,3}})]\wt g
\in M_{1\times(1+2\iota_3)}
\eqno(4.56)$$
(cf.~(4.24) and (4.33)), where $[\nu(x^{-\si})^T,\nu(\ol t'_{J_{1,3}})]\in
M_{(1+\iota_5)\times(1+2\iota_3)}$.
Note that if $J_0\ne\{0\}$, then by (4.5), $h_0=0$, and by (4.30)-(4.32),
$\th([t_0,\vsi_p])=0=[\th(t_0),\th(\vsi_p)]$ for
$p\in\wh I_3\cup I_{5,6}\cup \ol I_{2,6}$ if $\G_0\ne\{0\}$.
If $\G_0=\{0\}$, we have $\th([t_0,\vsi_p])=0=[\th(t_0),\th(\vsi_p)]$ for
$p\in I_4\cup J_{5,6}$ and
\par\noindent$
[\th(t_0),\th(\ol t_{J_{1,3}})]
\!=\![t'_0\!+\!\wt D^T\nu(x^{-\si})^T,\nu(\ol t'_{J_{1,3}})\wt g\!+\!E]
\!=\!-2E\!+\!\wt D^T[\nu(x^{-\si})^T,\nu(\ol t'_{J_{1,3}})]\wt g\!=\!0
\hfill(4.57)$\par\noindent
by (4.30) and (4.56).
For $p\in I_3$, we have
\begin{eqnarray*}\hspace{2cm}[\th(\ol t_{\{p,\ol p\}})^T,\th(\ol t_{\{p,\ol p\}})]&=&
b_p(g_p^{-1})^T[(\ol t'_{\{\nu(p),\ol{\nu(p)}\}})^T,\ol t'_{\{\nu(p),\ol{\nu(p)}\}}]b_pg_p^{-1}\\ &=&
b_p(g_p^{-1})^T\left(\begin{array}{cc}0&x^{\si_p}\\-x^{\si_p}&0\\ \end{array}\right)b_pg_p^{-1}\\
&=&\th([(\ol t_{\{p,\ol p\}})^T,\ol t_{\{p,\ol p\}}]).\hspace{5.2cm}(4.58)\end{eqnarray*}
By (4.31) and  (4.32), we also have
$$
[\th(t_{I_5})^T,\th(\vsi_{\ol I_5})]=B[(t'_{I_5})^T,(x^{-\si'})_{I_5}]B^{-1}
={\bf1}_{\ell_5}=\th([(t_{I_5})^T,\vsi_{\ol I_5}]).
\eqno(4.59)$$
Similarly, we can prove other identities.
Therefore, the last equation in (4.43) holds. This proves
Claim 1 in this case.
\par
{\it Case a.2}: Suppose $b_0\ne1$.
Multiplying $\tau$ by an element of ${\cal G}$ of the form in Case a.1,
we can suppose $\nu=\mbox{Id}_{\wh I}$ and $g'={\rm diag}(b_0,1,1,...,1,1)$.
Then it is straightforward to verify that (4.39) defines an isomorphism.
\par
{\it Case b}: Assume that $\nu(0)=p\ne0$. Then (4.9) shows that $J_0=\{0\}$
and $p\in I_1$.
By multiplying $\tau$ by an element of ${\cal G}$ of the form in Case a, we
can suppose $\nu(0)=p,\nu(p)=0$ and $\nu|_{I\bs\{p\}}=\mbox{Id}_{I\bs\{p\}}$
and $g'={\bf1}_{1+2\iota_6}$. In this case it is straightforward to verify
that
$$
\th: x^{\a,\vec i}\mapsto x^{\tau(\a)-\si_p,\vec i}
\eqno(4.60)$$
is an isomorphism $\AA\rar\AA'$.
\psp

``$\Longrightarrow$'' Suppose that $\theta: {\cal A}\rightarrow {\cal A}'$ is a Lie algebra isomorphism. We consider the following cases.
\par
{\it Case 1}. $\G_0\ne\{0\}$ and $\G'_0\ne\{0\}$.
\par
Note that (3.32) and (3.33) show that $\BB^{\rm F}/\BB^{\rm N}$ has dimension
$\iota_3$.
Thus we have
$$\iota_3=\iota'_3.
\eqno(4.61)$$
By the definition in (2.13), we have
$$
\si'_p=\si_p\;\;\;\;\for\;\;\;\;p\in\ol{0,\iota_3}.\eqno(4.62)$$
By Lemma 3.3, there exists a bijection $\tau_1:\pi(\G)\rar
\pi'(\G')$ such that
$$
\th(\BB_\mu)=\BB'_{\tau_1(\mu)}\;\;\; \for\ \mu\in\pi(\G).
\eqno(4.63)$$
Denote
$$\G_{0,3}=\{\a\in\G\,|\,(\pi(\a))_{\wh I_{1,3}}=0\}.
\eqno(4.64)$$
By (4.63) and Lemma 3.5, there exists a bijection
$\tau:\G\bs\G_{0,3}\rar
\G'\bs\G'_{0,3}$ such that
$$\th(x^{\a})=c_\a x^{\tau(\a)}\ \for\ \a\in\G\bs\G_{0,3}
\mbox{ and some }c_\a\in\mbb{F}\bs\{0\}.
\eqno(4.65)$$
\par
{\bf Claim 2}.
There exists a bijection $\nu:\wh I_{1,3}\rar\wh I'_{1,3}$
such that
$$
\th(x^{-\si_p})=d_px^{-\si_{\nu(p)}}\mbox{ for }p\in\wh I_{1,3}
\mbox{ and some }d_p\in\mbb{F}\bs\{0\}.
\eqno(4.66)$$
\par
By (3.10)-(3.11) and Lemma 3.4, we have
$$
\{u\in\BB^{\rm F}\,|\,[u,A_1]=0\}=
{\rm Span}\{x^{-\si_p}\,|\,p\in\wh I_{1,3}\}=
\{u\in\BB^{\rm F}\,|\,[u,A_2]=0\}.
\eqno(4.67)$$
Thus by Lemma 3.2,
$$
\{u\in\BB^{\rm F}\,|\,[u,\AA^{\rm N}]=0\}=
{\rm Span}\{x^{-\si_p}\,|\,p\in\wh I_{1,3}\}.
\eqno(4.68)$$
If $\iota_3=0$, by (4.61) and (4.68), there is nothing to prove.
Assume that $\iota_3>0$.
Suppose $p\in I_{1,3}$. By (2.15), $e_{[0]}\in\G$ for some $e\in\mbb{F}\bs\{0\}$.
By (3.36),
$$[x^{e_{[0]}-\si_p},x^{-e_{[0]}}]=-2ex^{-\si_p}.
\eqno(4.69)$$
Noting that $e_{[0]}-\si_p,-e_{[0]}\notin\G_{0,3}$, by (4.66),
$\th(x^{e_{[0]}-\si_p})\in\mbb{F} x^\a,\;\th(x^{-e_{[0]}})\in\mbb{F} x^\b$
for some $\a,\b\in\G'\setminus \G'_{0,3}$ (cf.~(4.64)).
We have
$$[x^\a,x^\b]=\sum_{q\in I'_{1,3}}(\a_q\b_{\ol q}-\a_{\ol q}\b_q)
x^{\si_q+\a+\b}+((2-\vt(\a,0))\b_0-\a_0(2-\vt(\b,0)))x^{\a+\b}
\eqno(4.70)$$
(cf.~(3.36)).
By (4.68),
$\th(x^{-\si_p})\in\sum_{q\in\wh I_{1,3}}\mbb{F} x^{-\si_q}$.
We claim that $\th(x^{-\si_p})\in\cup_{q\in\wh I_{1,3}}\mbb{F} x^{-\si_q}$.
If not, then by (4.69) and (4.70),  there exist $q,r\in\wh I'_{1,3}$ with $q\ne r$
such that $\si_q+\a+\b=-\si_r$, that is, $\b=-\a-\si_q-\si_r$.  If $q,r\ne0$, then (4.70) becomes
\begin{eqnarray*}[x^\a,x^\b]&=&(\a_q(1-\a_{\ol q})-\a_{\ol q}(1-\a_q))x^{-\si_r}+(\a_r(1-\a_{\ol r})-\a_{\ol r}(1-\a_r))x^{-\si_q}\\ &
=&(\a_q-\a_{\ol q})x^{-\si_r}+(\a_r-\a_{\ol r})x^{-\si_q},\hspace{7cm}(4.71)\end{eqnarray*}
with $\a_q-\a_{\ol q},\a_r-\a_{\ol r}\ne0$. We obtain that
$[x^\a,[x^\a,x^\b]]=2(\a_q-\a_{\ol q})(\a_r-\a_{\ol r})x^\a
\ne0$, which contradicts the fact that $[x^{e_{[0]}-\si_p},
[x^{e_{[0]}-\si_p},x^{-e_{[0]}}]]=0$. If $q\ne0, r=0$, then
(4.70) becomes
$$[x^\a,x^\b]
=(\a_q-\a_{\ol q})-2\a_0x^{-\si_q},
\eqno(4.72)$$
with $\a_q-\a_{\ol q},\a_0\ne0$ (cf.~(3.6)). So we have
$[x^\a,[x^\a,x^\b]]=-4\a_0(\a_q-\a_{\ol q})x^\a
\ne0$, which again contradicts the fact that $[x^{e_{[0]}-\si_p},
[x^{e_{[0]}-\si_p},x^{-e_{[0]}}]]=0$.
Thus $\th(x^{-\si_p})\in\cup_{q\in\wh I_{1,3}}\mbb{F} x^{-\si_q}$.
Now consider $\th(1_{\cal A})$. By (2.15),
$a_{[1]}\in\G$ for some $a\in\mbb{F}\bs\{0\}$. Then
$$[x^{a_{[1]}-\si_1},x^{-a_{[1]}}]=a1_{\cal A}.
\eqno(4.73)$$
As proved above, we have
$\th(1_{\cal A})\in\cup_{q\in\wh I_{1,3}}\mbb{F} x^{-\si_q}$. This proves
the claim.
\psp

Now we consider the following subcases.
\par
{\it Case 1(i)}: $\nu(0)=0$.
\par
Then $J_0=J'_0$ since $1_{\cal A}$  or $1_{{\cal A}'}$ is {\it ad}-semisimple
if and only if $J_0=\{0\}$ or $J'_0=\{0\}$ respectively.
For $p\in\wh I_{1,3}$, by (2.15) and (2.16),
we can fix $e_p\in\mbb{F}\bs\{0\}$ such that
$\eta_p=(e_p)_{[p]}\in\G$. Then $\eta_p\notin\G_{0,3}$.
\par
{\bf Claim 3}.
$$(\tau(\eta_p))_{\nu(q)}=0\qquad \for\ p\in\wh I_{1,3},\,q\in\wh J_{1,3},\,
q\ne p,\ol p.
\eqno(4.74)$$
\par
For $p,q\in I_{1,3}$, applying $\th$ to $[x^{-\si_p},x^{\eta_q}]=
-\d_{p,q}e_qx^{\eta_q}$, we obtain
$$d_p((\tau(\eta_q))_{\nu(p)}-(\tau(\eta_q))_{\ol{\nu(p)}})=-\d_{p,q}e_p\qquad\mbox{for }p,q\in I_{1,3}
\eqno(4.75)$$
 by (4.65) and (4.66), Applying $\th$ to $[x^{\eta_p},x^{\eta_q}]=0$ and using (4.75), we obtain
$$0=(\tau(\eta_p))_{\nu(q)}(\tau(\eta_q))_{\ol{\nu(q)}}-
(\tau(\eta_p))_{\ol{\nu(q)}}(\tau(\eta_q))_{\nu(q)}
=-(\tau(\eta_p))_{\nu(q)}e_q,
\eqno(4.76)$$
for $p,q\in I_{1,3},\,p\ne q.$
The above two equations show that (4.74) holds for $p,q\in I_{1,3}$
with $q\ne p,\ol p$. In the same way, from $[1_{\cal A},x^{\eta_p}]=0$, we obtain $(\tau(\eta_p))_0=0$.
Thus (4.74) holds for $p\in I_{1,3}$. Similarly, we can prove that
(4.74) holds for $p=0$.
\par
{\bf Claim 4}.
$\tau$ can be uniquely extended to a group isomorphism
$\tau:\G\rar\G'$ such that  $\tau(\si_p)=\si_{\nu(p)}$
for $p\in \wh{I}_{1,3}$.
\par
For any $\a\in\G$ with $\a,\a+\eta_0\notin\G_{0,3}$ (cf.~(4.64)),
we have
\begin{eqnarray*}& &c_{\eta_0+\a}
(2\a_0-e_0(2-\vt(\a,0)))x^{\tau(\eta_0+\a)}
\\ &=&\th([x^{\eta_0},x^\a])
=c_{\eta_0}c_\a[x^{\tau(\eta_0)},x^{\tau(\a)}]\\ &=&c_{\eta_0}c_\a(
(2-\vt(\tau(\eta_p),0))(\tau(\a))_0-(\tau(\eta_0))_0(2-\vt(\tau(\eta_0),0)))
x^{\tau(\eta_0)+\tau(\a)},\hspace{1.6cm}(4.77)\end{eqnarray*}
which implies
$$
\tau(\eta_0+\a)=\tau(\eta_0)+\tau(\a),
\eqno(4.78)$$
if $\a$ satisfies
$$
\a,\a+\eta_0\notin\G_{0,3}\mbox{ and }
2\a_0-e_0(2-\vt(\a,0))\ne0.
\eqno(4.79)$$
\par
For $\a,\b\in\G$ such that $\a,\b,\a+\b\notin\G_{0,3}$,
we have
\begin{eqnarray*}& &
\sum_{p\in I_{1,3}}(\a_p\b_{\ol p}-\a_{\ol p}\b_p)c_{\si_p+\a+\b}
x^{\tau(\si_p+\a+\b)}\\ & &
+((2-\vt(\a,0))\b_0-\a_0(2-\vt(\b,0)))c_{\a+\b}x^{\tau(\a+\b)}
\\ &=&c_\a c_\b(
\sum_{p\in I'_{1,3}}((\tau(\a))_{\nu(p)}(\tau(\b))_{\ol{\nu(p)}}-
(\tau(\a))_{\ol{\nu(p)}}(\tau(\b))_{\nu(p)})
x^{\si_{\nu(p)}+\tau(\a)+\tau(\b)}\\ & &+((2-\vt(\tau(\a),0))(\tau(\b))_0-(\tau(\a))_0(2-\vt(\tau(\b),0)))x^{\tau(\a)+\tau(\b)})\hspace{2.7cm}(4.80)\end{eqnarray*}
by applying $\th$ to (3.36). Assume that all $\a,\a+\eta_0,2\a,2\a+\eta_0$ satisfy (4.79).
Taking $\b=\eta_0+\a$ in (4.80) and using (4.78), we see that all terms
in (4.80) vanish except the last terms in the both sides. Thus we obtain
$$
\tau(2\a)=2\tau(\a),
\eqno(4.81)$$
if (4.79) holds for $\a,\a+\eta_0,2\a,2\a+\eta_0$ and
$2-\vt(\a,0)\ne0.$
Since we can replace $\eta_0$ by $k\eta_0$ for any $0\ne k\in\mbb{Z}$,  (4.81)
holds for all $\a\in\G\bs\G_{0,3}$ with $2-\vt(\a,0)\ne0.$
\par
Denote
$$\Si=\{\a\in\G\bs\G_{0,3}\,|\,2-\vt(\a,0)\ne0\}.\eqno(4.82)$$
Now take $\a,\b\in\G$ such that
$$\a,\b,\a+\b\in \Si\mbox{ and }
(2-\vt(\a,0))\b_0-\a_0(2-\vt(\b,0))\ne0.
\eqno(4.83)$$
Then (4.80) implies that
$$
\begin{array}{ll}\dis
\tau(\a+\b)=\tau(\a)+\tau(\b)+\sum_{p\in I_{1,3}}k^{(p)}_{\a,\b}\si_p,
\mbox{ where}
\vs{4pt}\\ \dis
k^{(p)}_{\a,\b}=0,1\mbox{ such that }
\sum_{p\in I_{1,3}}k^{(p)}_{\a,\b}=0,1.
\end{array}
\eqno(4.84)$$
We claim that $\tau(\a+\b)=\tau(\a)+\tau(\b)$ if $\a,\b\in \Si$ and if
the pairs $(\a,\b),(2\a,2\b)$ satisfy (4.83).
Assume that $k_{\a,\b}^{(q)}=1$ for some $q\in I_{1,3}$.
Then we obtain
\begin{eqnarray*}\tau(2\a)+\tau(2\b)+\sum_{p\in I_{1,3}}k_{2\a,2\b}^{(p)}\si_p
&=&\tau(2\a+2\b)=\tau(2(\a+\b))=2\tau(\a+\b)\\ &=&2(\tau(\a)+\tau(\b)+\sum_{p\in I_{1,3}}k_{\a,\b}^{(p)}\si_p),\hspace{3.2cm}(4.85)\end{eqnarray*}
from this we obtain $k_{2\a,2\b}^{(q)}=2k_{\a,\b}^{(q)}>1$, which is a contradiction
to (4.84).
\par
For any $\a,\b\in \Si$, we can always choose $\g\in \Si$ such that
the pairs
$$(\a+\b,\g),\;(2\a+2\b,2\g),\;(\a,\b+\g),\;(2\a,2\b+2\g),\;
(\b,\g),\;(2\b,2\g)\eqno(4.86)$$
satisfy (4.83). Hence
$$\tau(\a+\b)+\tau(\g)=\tau(\a+\b+\g)=\tau(\a)+\tau(\b+\g)=
\tau(\a)+\tau(\b)+\tau(\g),\eqno(4.87)$$
which  gives rise to
$$\tau(\a+\b)=\tau(\a)+\tau(\b).\eqno(4.88)$$
In this way, we can also prove $\tau(\a+\b)=\tau(\a)+\tau(\b)$ for all
$\a,\b\in\G\bs\G_{0,3}$ such that $\a+\b\in\G\bs\G_{0,3}$.
Since any element $\g\in\G$ can be written as $\g=\a+\b$ for some
$\a,\b\in\G\bs\G_{0,3}$, this shows that $\tau$ can be uniquely extended
to a group isomorphism $\tau:\G\rar\G'$ such that $\tau(\si_p)=\si_{\nu(p)}$
for $p\in \wh{I}_{1,3}$.
\par
{\bf Claim 5}.
There exists $g={\rm diag}(b_0,g_1,...,g_{\iota_3})\in GL_{1+2\iota_3}$,
where
$$b_0\in\mbb{F}\bs\{0\},\;\;
g_p=\left(\begin{array}{cc}1-a_p&1-a_p-b_p\\ a_p&a_p+b_p\end{array}\right)
\in GL_2,
\eqno(4.89)$$
such that
$$
\nu(\tau(\a))_{\wh J_{1,3}}=\a_{\wh I_{1,3}}g\mbox{ and }
\vt(\tau(\a),0)=\vt(\a,0)\qquad \for\ \a\in\G,
\eqno(4.90)$$
(cf.~(4.10)).
\par
Using the fact that $\tau$ is a group isomorphism and applying $\th$ to
$${[1_{\cal A},x^\a]=2\a_0x^\a},\;\;[x^{-\si_p},x^\a]=(\a_{\ol p}-\a_p)x^\a,\eqno(4.91)$$
$${[x^\a,x^{-\a}]=-4\a_0},\;\;[x^\a,x^{-\a-\si_p}]=(\a_p-\a_{\ol p})-2\a_0x^{-\si_p}\eqno(4.92)$$
for $\a\in\G\bs\G_{0,3}$ and $p\in I_{1,3},$
we obtain
$$2\a_0=2(\tau(\a))_0d_0,\;\;\a_{\ol p}-\a_p=d_p((\tau(\a))_{\ol{\nu(p)}}-(\tau(\a))_{\nu(p)}),\eqno(4.93)$$
$$-4\a_0d_0=-4(\tau(\a))_0c_\a c_{-\a},\;\;-2\a_0d_p=-2(\tau(\a))_0c_\a c_{-\a-\si_p},\eqno(4.94)$$
$$(\a_p-\a_{\ol p})d_0=c_\a c_{-\a-\si_p}((\tau(\a))_{\nu(p)}-(\tau(\a))_{\ol{\nu(p)}})\eqno(4.95)$$
for $\a\in\G\bs\G_{0,3}$ and $p\in I_{1,3}.$
Comparing the coefficients of $x^{\si_{\nu(p)}+\tau(\a)+\tau(\b)}$
in (4.80), we obtain
$$(\a_p\b_{\ol p}-\a_{\ol p}\b_p)c_{\si_p+\a+\b}
=c_\a c_\b((\tau(\a))_{\nu(p)}(\tau(\b))_{\ol{\nu(p)}}
-(\tau(\a))_{\ol{\nu(p)}}(\tau(\b))_{\nu(p)}),\eqno(4.96)$$
\begin{eqnarray*}& &((2-\vt(\a,0))\b_0-\a_0(2-\vt(\b,0)))c_{\a+\b}\\
&=&c_\a c_\b((2-\vt(\tau(\a),0))(\tau(\b))_0-(\tau(\a))_0(2-\vt(\tau(\b),0))),\hspace{4cm}(4.97)\end{eqnarray*}
\begin{eqnarray*}\qquad& &((1-\a_p)(1-\b_{\ol p})-(1-\a_{\ol p})(1-\b_p))c_{-\si_p-\a-\b}
\\ &=&c_{-\a-\si_p} c_{-\b-\si_p}
((1-(\tau(\a))_{\nu(p)})(1-(\tau(\b))_{\ol{\nu(p)}})\\ & &
-(1-(\tau(\a))_{\ol{\nu(p)}})(1-(\tau(\b))_{\nu(p)})),\hspace{7cm}(4.98)\end{eqnarray*}
\begin{eqnarray*}& &(-\vt(\a,0)\b_0+\a_0\vt(\b,0))c_{-\a-\b-2\si_p}\\
&=&c_{-\a-\si_p} c_{-\b-\si_p}(-\vt(\tau(\a),0)(\tau(\b))_0
+(\tau(\a))_0\vt(\tau(\b),0))\hspace{4.3cm}(4.99)\end{eqnarray*}
for $\a,\b,\a+\b\in\G\bs\G_{0,3}$,
where (4.98) and (4.99) are obtained from (4.96) and (4.97) through replacing
$\a$ and $\b$ by $-\a-\si_p$ and $-\b-\si_p$, respectively.
Set $\b=\eta_p$ in (4.96). By (4.74), (4.75) and  the second equation in (4.93), we obtain
\begin{eqnarray*}-\a_p e_pc_{\si_p+\eta_p+\a}
&=&c_{\eta_p}c_\a((\tau(\a))_{\nu(p)}(-d_p^{-1}e_p+(\tau(\eta_p))_{\nu(p)})
\\ &&-(d_p^{-1}(\a_{\ol p}-\a_p)+(\tau(\a))_{\nu(p)})(\tau(\eta_p))_{\nu(p)})
\\  &=&-d_p^{-1}c_{\eta_p}c_\a(e_p(\tau(\a))_{\nu(p)}+(\a_{\ol p}-\a_p)(\tau(\eta_p))_{\nu(p)}),
\hspace{2.8cm}(4.100)\end{eqnarray*}
if $\a,\a+\eta_p\in\G\bs\G_{0,3}.$
Set $\b=\eta_p$ in (4.99) and note that
$\vt(\eta_p,0)=e_p$. We obtain
$$
\a_0e_pc_{-2\si_p-\eta_p-\a}=c_{-\si_p-\eta_p}
c_{-\si_p-\a}(\tau(\a))_0\vt(\tau(\eta_p),0)
\mbox{ if }\a,\a+\eta_p\in\G\bs\G_{0,3}.
\eqno(4.101)$$
Assume that $\a_0\ne0$, which implies $\a,\a+\eta_p\in\G\bs\G_{0,3}$.
By the first equation in (4.93) and the second equation in (4.94),
we have $c_{\si_p+\eta_p+\a}c_{-2\si_p-\eta_p-\a}=c_\a c_{-\si_p-\a}
=d_0d_p$. Since $\eta_p$ is fixed for each $p$,
 (4.100), (4.101), the second equation in (4.93) and the fact
that $\tau(\si_p)=\si_{\nu(p)}$ enable us to  deduce
$$((\tau(\a))_{\nu(p)},(\tau(\a))_{\ol{\nu(p)}})=(\a_p,\a_{\ol p})\left(\begin{array}{cc}1-a_p&1-a_p-b_p\\ a_p&a_p+b_p\\ \end{array}\right)\eqno(4.102)$$
 for $p\in I_{1,3}$ and  some $a_p,b_p\in\mbb{F}$ with $b_p\ne0$, where the matrix is denoted as $g_p$ in (4.18).
Since $\tau$ is a group isomorphism, (4.102) must hold for all $\a\in\G$.
Note that the determinant $|g_p|=b_p\ne0$. We also have
$$
(\tau(a))_0=b_0\a_0,\mbox{ where }b_0=d_0^{-1}\ne0,
\eqno(4.103)$$
by the first equation in (4.93). This proves the first equation in (4.90).
\par

Now we have
$$
\begin{array}{ll}
d_p|g_p|=1,&
\vs{4pt}\\
c_\a c_{-\a}=d_0^2&\mbox{ if }\a_0\ne0,
\vs{4pt}\\
c_a c_{-\si_p-\a}=d_0d_p&\mbox{ if }\a_0\ne0,
\vs{4pt}\\
c_{\si_p+\a+\b}=c_\a c_\b|g_p|
&\mbox{ if }\a_p\b_{\ol p}-\a_{\ol p}\b_p\ne0,
\\
\end{array}
\eqno(4.104)$$
by (4.93), (4.94) and (4.96), where $\a,\b,\a+\b\in\G\bs\G_{0,3}$.
{}From this one can deduce that
$\chi:\G\bs\G_{0,3}\rar\mbb{F}^{\times}$ defined by $\chi(\a)=b_0c_\a$ can
be uniquely extended to a multiplicative function $\chi:\G\rar\mbb{F}^{\times}$.
Then by (4.103) and (4.97), we obtain
the second equation of (4.90).
\par

{\bf Claim 6}.
We claim that $\nu(I_k)=I'_k$ for $k=1,2,3$.
\par
So suppose $\iota_3\ge1$.
Take
$${\cal M}={\rm Span}\{x^\a\,|\,\a\in\G\bs\G_{0,3},\,\a_0=0\},\eqno(4.105)$$
$$
\begin{array}{ll}
{\cal M}_1\!\!\!\!&=\{u\in\AA\,|\,[u,{\cal M}]\subset\BB\}
\vs{4pt}\\ &
=\BB+{\rm Span}\{x^{\a,\vec i},x^{\b,\vec j}\,|\,
\a=\a_{I_{4,5}},\,|\vec i|=1,\b=\b_{I_{4,5}},\vec j=\vec j_{I_5\cup J_6}\},
\end{array}
\eqno(4.106)$$
$$
\begin{array}{ll}
{\cal M}_2\!\!\!\!&=
\BB+\{u\in{\cal M}_1\,|\,[x^{-\si_p},u]=0\mbox{ for }p\in\wh I_{1,3}\}
\vs{4pt}\\ &
=\BB+{\rm Span}\{x^{\a,1_{[q]}},x^{\b,\vec j}\,|\,\a=\a_{I_{4,5}},
q\in\ol I_{4,5},\b=\b_{I_{4,5}},\vec j=\vec j_{I_5\cup J_6}\},
\end{array}
\eqno(4.107)$$
$${\cal M}_1^{(p)}=\BB+\{u\in{\cal M}_1\,|\,[x^{-\si_p},u]=0\}\mbox{ for }p\in\wh I_{1,3}.
\eqno(4.108)$$
Then ${\cal M}_2$ is a Lie algebra and ${\cal M}_1$ is a ${\cal M}_2$-module
such that ${\cal M}_1^{(p)}$ is a submodule for $p\in I_{1,3}$.
Note that the quotient module ${\cal M}_1/{\cal M}_1^{(p)}$ is zero if
$p\in I_1$, is a cyclic ${\cal M}_2$-module (with generator $t_{\ol p}$)
if $p\in I_2$, and is not cyclic (with two generators $t_p,t_{\ol p}$) if
$p\in I_3$. Applying $\th$ to the above sets, we obtain
$$
(\nu(I_1),\nu(I_2),\nu(I_3))=(I'_1,I'_2,I'_3)
\mbox{ and so }(\ell_1,\ell_2,\ell_3)=(\ell'_1,\ell'_2,\ell'_3).
\eqno(4.109)$$
\par
{\bf Claim 7}. In (4.102), $a_p=0$ if $p\in I_2$.
\par
Assume that $p\in I_2$.
Write
$$\th(t_{\ol p})=b t'_{\nu(p)}+
\sum_{(0,1_{[\ol{\nu(p)}]})\ne(\b,\vec j)\in\G'\times\JJ'}b_{\b,\vec j}
x^{\b,\vec j}\mbox{ for some }b,b_{\b,\vec j}\in\mbb{F}.
\eqno(4.110)$$
We have
$$
\a_pc_{\si_p+\a}x^{\tau(\a)+\si_{\nu(p)}}=\th([x^\a,t_{[\ol p]}])
=[\th(x^\a),\th(t_{[\ol p]})]
=c_a b(\tau(\a))_{\nu(p)}x^{\tau(\a)+\si_{\nu(p)}}+...,
\eqno(4.111)$$
where missed terms do not contain $x^{\tau(\a)+\si_{\nu(p)}}$.
Thus by (4.102), we obtain
$$
\a_pc_{\si_p+\a}=c_a b((1-a_p)\a_p+a_p\a_{\ol p}).
\eqno(4.112)$$
By (2.15), we can choose $\a\in(\mbb{F} 1_{[\ol p]}\cap\G)\bs\{0\}$, then (4.112)
proves the claim.
\par
This proves (4.1)-(4.3).
\par
By (3.9), (3.11) and Lemma 3.2, we have
$$
\begin{array}{lll}
\th(x^{-\si_p,1_{[\ol p]}})\in\AA'{}^{\rm F}&\subset&{\rm
Span}(A'_0\cup A'_2) \vs{4pt}\\ &=&\dis\sum_{q\in
I'_{1,3}}\mbb{F}x^{-\si_p}+\sum_{r\in I'_{4,5}}
\mbb{F}x^{-\si_p,1_{[\ol p]}}+{\rm Span}(A'_2),
\\ \end{array}
\eqno(4.113)$$
for $p\in I_{4,5}$.
Thus,
$$\th((x^{-\si})_{I_{4,5}})\equiv
(\nu(x^{-\si'})_{\wh{I'}_{1,3}})\wt h+(x^{-\si'})_{I'_{4,5}}\wt f
\,\ ({\sc\,}{\rm mod\ \,}{\rm Span}(A_2'){\sc\,})
\eqno(4.114)$$
(cf.~(4.27)) for some
$$\wt h=(a_{p,q})\in M_{(1+\iota'_3)\times(\ell_4+\ell_5)}
\mbox{ and }\wt f =(b_{p,q})\in GL_{\ell_4+\ell_5}.
\eqno(4.115)$$
\par
{\bf Claim 8}.
$$a_{p,q}=0\mbox{ if }p\in \wh{I'}_{2,3}\mbox{ or }p=0,J'_0\ne\{0\},
\eqno(4.116)$$
$$b_{p,q}=0\mbox{ for }p\in I'_5,\,q\in I_4.
\eqno(4.117)$$
$$
(\ell_4,\ell_5)=(\ell'_4,\ell'_5),\ \ \ \si=\si'.
\eqno(4.118)$$
\par
Recall (3.35). Denote
$$C_\AA(C(\BB))=\{u\in\AA\,|\,[u,C(\BB)]=0\},
\eqno(4.119)$$
the centralizer of $C(\BB)$.
By (3.35), we note that $x^{-2\si_p}\in C(\BB)$ for $p\in I_{4,5}$.
It is straightforward to verify
$$\{t_p\,|\,p\in I_3\cup\ol I_{2,3}\}\subset
C_\AA(C(\BB))\subset{\rm Span}
\{x^{\a,\vec i}\,|\,\vec i_{\ol I_{4,5}}=0\},
\eqno(4.120)$$
and $t_0\in C_\AA(C(\BB))$ if $J_0\ne\{0\}$.
For $p\in\wh I_{1,3}$,  (4.120) implies that
${\rm ad}_{x^{-\si_p}}|_{C_\AA(C(\BB))}$ is
semisimple if and only if $p\in I_1$ or $p=0,J_0=\{0\}$,
and ${\rm ad}_{x^{-\si_q,1_{[\ol q]}}}|_{C_\AA(C(\BB))}$ is
semisimple for $q\in I_{4,5}$ by (4.120). Moreover, by (3.1),
for $q\in I_{4,5}$,
${\rm ad}_{x^{-\si_q,1_{[\ol q]}}}$ is semisimple if and only if
$q\in I_4$. We obtain (4.116) and (4.117). In particular,
(4.117) and the definition of $\si$ in (2.13) imply (4.118).
\par
For any $\a\in\G$, we denote
$$
\ol\a=(2\a_0,\a_{\ol1}-\a_1,...,\a_{\ol\iota_3}-\a_{\iota_3})
=\a_{\wh J_{1,3}}P\in\mbb{F}^{1+\iota_3}\eqno(4.121)$$
(cf.~(4.15)). For $\a\in\G\bs\G_{0,3}$, applying $\th$ to
$$
-\a_{I_{4,5}} x^\a=[(x^{-\si})_{I_{4,5}},x^\a],
\eqno(4.122)$$
(cf.~(3.17) and (4.27)),  using (4.65), (4.66) and (4.114), and noting that
$[A'_2,\BB']=0$, we obtain
\begin{eqnarray*}\qquad-\a_{I_{4,5}}c_\a x^{\tau(\a)}&=&
c_\a[\nu((x^{-\si})_{\wh{I'}_{1,3}})\wt h+(x^{-\si})_{I'_{4,5}}\wt f,
x^{\tau(\a)}]\\ &=&c_\a
(\nu(\ol{\tau(\a)})_{\wh{I'}_{1,3}}\wt h-(\tau(\a))_{I_{4,5}}\wt f)x^{\tau(\a)}\\
&=&c_\a(\nu((\tau(\a))_{\wh J_{1,3}})P\wt h-(\tau(\a))_{I_{4,5}}\wt f)
x^{\tau(\a)}\hspace{3.7cm}(4.123)\end{eqnarray*}
(cf.~(4.121)). This gives
$$(\tau(\a))_{I_{4,5}}=\a_{I_{4,5}}f+\a_{\wh J_{1,3}}h\mbox{ with }f=\wt f^{-1},\,
h=gP\wt h\wt f^{-1},
\eqno(4.124)$$
for all $\a\in\G\bs\G_{0,3}$ and thus for all $\a\in\G$ since $\tau$ is an
isomorphism, where
$$
g={\rm diag}(b_0,g_1,...,g_{\iota_3})\in GL_{1+2\iota_3}
\eqno(4.125)$$
({cf.~(4.102) and (4.103)).
Clearly $f$ and $h$ are of the forms in (4.4) and (4.6), respectively.
Assume that $\ell_6=\ell'_6$, then
this proves that $\tau$ has the form $g'_\nu$ in (4.12), and
(4.90) implies that $g'\in{\cal G}''$ (cf.~(4.7)), and thus
the theorem is proved in this case. It remains to prove
\par
{\bf Claim 9}. $\ell_6=\ell'_6$.
\par
Recall (4.105) for the definition of ${\cal M}$. By (3.1) we see that the
centralizer of ${\cal M}$ in $\AA$ is
$$
{\cal C}=C_\AA({\cal M})={\rm Span}\{x^{\a,\vec i}\,|\,\a=\a_{I_{4,5}},\,
\vec i=\vec i_{I_5\cup J_6}\}.
\eqno(4.126)$$
The center of ${\cal C}$ is
$$
{\cal D}=C({\cal C})=
{\rm Span}\{x^{\a,\vec i}\,|\,\a=\a_{I_{4,5}},\vec i=\vec i_{I_5}\}
\mbox{ which is a domain ring}.
\eqno(4.127)$$
By exchanging positions of $\AA$ and $\AA'$ if necessary, we can suppose
$\ell_6\le\ell'_6$.
Note that in ${\cal C}$,  the formula (1.1) holds.
By the proof of sufficiency, there exists an embedding $\ol\th:\AA\rar\AA'$
such that
$$
\ol\th(x^\a)=c_\a x^\a,\;\;
\ol\th((x^{-\si})_{I_{4,5}})=
(\nu(x^{-\si'})_{\wh{I'}_{1,3}})\wt h+(x^{-\si'})_{I'_{4,5}}\wt f\;\;\;
\for\;\;\a\in\G
\eqno(4.128)$$
(cf.~(4.31), (4.39), (4.114), and noting that the proof of Claim 5
shows that $\chi(\a)=b_0 c_\a$ is a multiplicative function). Thus
by identifying $\AA$ with $\ol\th(\AA)$, we can assume that $\AA$ is
a subalgebra of $\AA'$ such that there exists isomorphism $\th$ satisfying
$$
\th(x^\a)=x^\a,\;\;\th(x^{-\si_p,1_{[\ol p]}})\equiv x^{-\si_p,1_{[\ol p]}}
\,\ ({\sc\,}{\rm mod\ \,}{\rm Span}(A_2'){\sc\,})
\;\;\;\for\;\;\;\a\in\G,p\in I_{4,5}.
\eqno(4.129)$$
By restricting $\th$ to ${\cal C}$, we want to prove
$$
\th(t_p)=t_p+c_p,
\;\;\;\for\;\;\;p\in I_5\mbox{ and some }c_p\in\mbb{F},
\eqno(4.130)$$
$$
\th(x^{\a,\vec i} t^{\vec j})=
x^\a\prod_{p\in I_5}(\th(t_p))^{i_p}\prod_{q\in J_6}(\th(t_q))^{j_q}
\;\;\;\for\;\;\;\a=\a_{I_{4,5}},\vec i=\vec i_{I_5},\vec j=\vec j_{J_6},
\eqno(4.131)$$
To prove (4.130),
first by (4.127), we have $c_p=\th(t_p)-t_p\in {\cal D}'$. Then by
(4.129), we have
$$
[x^{-\si_q,1_{[\ol q]}},c_p]=
\th([\th^{-1}(x^{-\si_q,1_{[\ol q]}}),t_p])-[x^{-\si_q,1_{[\ol q]}},t_p]
=0,
\eqno(4.132)$$
from which, we obtain that $c_p\in\mbb{F}$. Thus we have (4.130).
Similarly, we have
$$
\th(x^{\a,1_{[p]}})=x^\a(t_p+c_{\a,p})
\;\;\;\for\;\;\;p\in I_5\mbox{ and some }c_{\a,p}\in\mbb{F}.
\eqno(4.133)$$
By considering $[x^\a,t_pt_{\ol p}]$, we see that $c_{\a,p}=c_p$, and
obtain
$$
\th(t_pt_{\ol p})=(t_p+c_p)t_{\ol p}+u_p
\;\;\;\for\;\;\;p\in I_5
\mbox{ and some }u_p\in C_{\AA'}({\cal D}').
\eqno(4.134)$$
{}From this and (4.130), we can deduce
$$
\th(x^{\a,1_{[p]}})=x^\a(t_p+c_p)\;\;\;\for\;\;\;p\in I_5.
\eqno(4.135)$$
Similar to (4.134), we have
$$
\th(x^{-\si_p,1_{[p]}+1_{[\ol p]}})=
x^{-\si_p,1_{[\ol p]}}(t_p+c_p)+u'_p
\;\;\;\for\;\;\;p\in I_5
\mbox{ and some }u'_p\in C_{\AA'}({\cal D}').
\eqno(4.136)$$
Now from  (4.130), (4.133)-(4.136), we can obtain
(4.131) by induction on $|\vec i|$ in case $\vec j=0$.
\par

Assume that (4.131) holds for all $\vec j$ with $|\vec j|<n$,
where $n\ge1$.
We denote by $X_{\a,\vec i,\vec j}$ 
the difference between the left-hand side and the right-hand 
side of (4.131). 
Then the inductive assumption says that $X_{\a,\vec i,\vec j}=0$  
if $|\vec j|<n$. Now suppose $|\vec j|=n$. Say $j_r\ge1$
for some $r\in I_6$ (the proof is similar if $r\in\ol I_6$). 
Let $\vec k=\vec j-1_{[r]}+1_{[\ol r]}$.
Then we have
\begin{eqnarray*}[\th(t_r),X_{\a,\vec i,\vec k}]
&=&\th([t_r,x^{\a,\vec i}t^{\vec k}])
-\th([t_r,\th^{-1}(x^\a)])
\prod_{p\in I_5}(\th(t_p))^{i_p}\prod_{q\in J_6}(\th(t_q))^{k_q}
\\ && -x^\a[\th(t_r),\prod_{p\in I_5}(\th(t_p))^{i_p}\prod_{q\in J_6}
(\th(t_q))^{k_q}]\\ &=&(j_{\ol r}+1)(
\th(x^{\a,\vec i}t^{\vec j-1_{[r]}})-x^\a\prod_{p\in I_5}(\th(t_p))^{i_p}\prod_{q\in J_6}(\th(t_q))^{j_q-\d_{q,r}})\\ &=&(j_{\ol r}+1)X_{\a,\vec i,\vec j-1_{[r]}}=0,\hspace{7.5cm}(4.137)\end{eqnarray*}
where the first equality follows from (1.1), and the second equality follows
from (1.1) and (4.129).
By (1.1) and (4.137), we obtain
$$
[\th(t^2_r),X_{\a,\vec i,\vec k}]=
\th([t^2_r,\th^{-1}(X_{\a,\vec i,\vec k})])
=2\th(t_r[t_r,\th^{-1}(X)])
=0.\eqno(4.138)$$
On the other hand, exactly similar to (4.137), we have
$$
[\th(t^2_r),X_{\a,\vec i,\vec k}]=
2(j_{\ol r}+1)X_{\a,\vec i,\vec j}.
\eqno(4.139)$$
Now (4.138) and (4.139) show that $X_{\a,\vec i,\vec j}=0$.
This proves (4.131). By (4.130), (4.131) and identifying
${\cal D}$ with ${\cal D}'$ using the isomorphism, we see that $\th$
is an associative algebra isomorphism ${\cal C}\to{\cal C}'$ over
the domain ring ${\cal D}$.
{}From this we obtain $\ell_6=\ell'_6$ since $2\ell_6$
is the transcendental degree of ${\cal C}$ over the domain ring
${\cal D}$.

\par
{\it Case 1(ii)}: Suppose $\nu(0)\ne0$.
\par
{\bf Claim 10}. $J'_0=\{0\}$.
\par
Take
$$\wt\G=({\bf1}_{1+2\iota_6})_{\nu^{-1}}(\G'),
\eqno(4.140)$$
and define a Lie algebra $\wt A={\cal K}(\vec\ell,\si,\wt G,\JJ)$.
Denote
$$
{\cal E}={\rm Span}\{x^{\a,\vec i}\in\BB\,|\,\vec i_{\wh J_{1,3}}=0\}.
\eqno(4.141)$$
Similarly, we have $\wt{\cal E}$ and ${\cal E}'$.
Clearly
$${\cal E}
\cong{\cal K}(\vec n,\si,\G,\JJ)\mbox{ with }
\vec n=(\iota_3,0,0,\ell_4,\ell_5,0).
\eqno(4.142)$$
Thus by the sufficiency, there exists
$$
\wt\th:{\cal E}'\cong\wt{\cal E}
\mbox{ which maps }
x^{-\si_{\nu(p)}}\mapsto \wt x^{-\si_p}\mbox{ for }p\in\wh I_{1,3}.
\eqno(4.143)$$
Also we have
$$
{\cal E}=\mbox{ the span of common eigenvectors of }\ad_{x^{-\si_p}}
\mbox{ for }p\in\wh I_{1,3}.
\eqno(4.144)$$
Thus
$$
\th'=\th|_{\cal E}:{\cal E}\cong{\cal E}'
\mbox{ maps }
x^{-\si_p}\mapsto x^{-\si_{\nu(p)}}\mbox{ for }p\in\wh I_{1,3}.
\eqno(4.145)$$
Therefore $\th''=\wt\th\cdot\th':{\cal E}\cong\wt{\cal E}$
which maps $x^{-\si_p}\mapsto x^{-\si_p}$ for $p\in\wh I_{1,3}$.
By the proof of Case 1(i), there exists $g'$ of the some suitable form
in (4.17) (with the corresponding sets $I_2,I_3$ being empty) such that
$\th''$ is determined by $\tau''=g'_{\rm Id}:\G\cong\wt\G$.
Thus $g'_\nu:\G\cong\G'$.
For simplicity, we assume $g'={\bf1}_{1+2\iota_6},\nu(0)=p_0,\nu(p_0)=0$
and $\nu|_{I\bs\{p_0\}}={\rm Id}|_{I\bs\{p_0\}}$ (the proof is similar
in general). For $\a\in\G$, we denote $\a'=\tau''(\a)$. Then
$$
\a'_0={1\over2}(\a_{\ol p_0}-\a_{p_0}),
\a'_{p_0}=-{1\over2}\vt(\a,0)-\a_0,\a'_{\ol p_0}=-{1\over2}\vt(\a,0)+\a_0,
\a'_p=\a_p,
\eqno(4.146)$$
for $p\ne0,p_0,\ol p_0,$
and
$$
\th(x^\a)=x^{\a'-\si_{p_0}}\mbox{ for }\a\in\G.
\eqno(4.147)$$
Suppose $J'_0\ne\{0\}$. Then we have $p_0\in I_2$ as in the proof of Claim 6. Write
$$
\th^{-1}(t'_0)=bx^{-\si_{p_0},1_{[\ol p_0]}}+
\sum_{(-\si_{p_0},1_{[\ol p_0]})\ne(\b,\vec j)\in\G\times\JJ',|\vec j|\le1}
b_{\b,\vec j}x^{\b,\vec j}.
\eqno(4.148)$$
If $b_{\b,\vec j}\ne0$ for some $\b\notin\mbb{F}\si_{p_0}$, then we can find some
$\a\in\G$ with $\a_{p_0}=\a_{\ol p_0}$ such that $[x^\a,x^{\b,\vec j}]\ne0$.
But by (4.146), $[t'_0,x^{\a'-\si_{p_0}}]=0$, which leads a contradiction.
Thus $\b\in\mbb{F}\si_{p_0}$.
Similarly, we have $\vec j=0$ if $b_{\b,\vec j}\neq 0$. Thus
(4.148) can be rewritten as
$$
\th^{-1}(t'_0)=bx^{-\si_{p_0},1_{[\ol p_0]}}+\sum_{i\in\sF}b_ix^{i\si_{p_0}}.
\eqno(4.149)$$
For any $\a\in\G$ with $\a_0=0$ and $\a_{p_0}\ne\a_{\ol p_0}$,
applying $\ad_{x^\a}$ to (4.149), we obtain
\begin{eqnarray*}& &
-2\a'_0\th^{-1}(x^{\a'-\si_{p_0},1_{[0]}})-\vt(\a',0)x^\a
\\ &=&\th^{-1}(-2\a'_0x^{\a'-\si_{p_0},1_{[0]}}-\vt(\a',0)x^{\a'-\si_{p_0}})
\\ &=&\th^{-1}([x^{\a'-\si_{p_0}},t'_0])=[x^\a,\th^{-1}(t'_0)]\\
&=&b((\a_{p_0}-\a_{\ol p_0})x^{\a,1_{[\ol p_0]}}+\a_{p_0}x^\a)
-\sum_{i\in\sF}ib_i(\a_{p_0}-\a_{\ol p_0})x^{\a+(i+1)\si_{p_0}}.
\hspace{3cm}(4.150)\end{eqnarray*}
Thus
\begin{eqnarray*}\qquad& &-2\a'_0\th^{-1}(x^{\a'-\si_{p_0},1_{[0]}})\\
&=&-2\a'_0bx^{\a,1_{[\ol p_0]}}
+(\a_{p_0}b+\vt(\a',0))x^\a
+2\a'_0\sum_{i\in\sF}ib_ix^{\a+(i+1)\si_{p_0}}.\hspace{2.5cm}(4.151)\end{eqnarray*}
Applying $\ad_{x^{-\a-2\si_{p_0}}}$ to (4.151) (note that $\tau''(\si_{p_0})
=-\si_{p_0}$ by (4,131)), we obtain
\begin{eqnarray*}& &-8{\a'}^2_0(bx^{-\si_{p_0},1_{[\ol p_0]}}+\sum_{i\in\mbb{F}}b_ix^{i\si_{p_0}})
-2\a'_0(4+\vt(\a',0))x^{-\si_{p_0}}\\ &=&-8{\a'}^2_0\th^{-1}(t'_0)-2\a'_0(4+\vt(\a',0))x^{-\si_{p_0}}\\ &=&-2\a'_0\th^{-1}(4\a'_0t'_0+(4+\vt(\a',0))\\ &=&
-2\a'_0\th^{-1}([x^{-\a'+\si_{p_0}},x^{\a'-\si_{p_0},1_{[0]}}])\\
&=&[x^{-\a-2\si_{p_0}},-2\a'_0\th^{-1}(x^{\a'-\si_{p_0},1_{[0]}})]
\\ &=&-2\a'_0b(-2(\a_{p_0}-\a_{\ol p_0})x^{-\si_{p_0},1_{[\ol p_0]}}
+(-\a_{p_0}+2)x^{-\si_{p_0}})\\ & &-2(\a_{p_0}b+\vt(\a',0))(\a_{p_0}-\a_{\ol p_0})x^{-\si_{p_0}}
+2\a'_0\sum_{i\in\mbb{F}}i(i-1)b_i(\a_{p_0}-\a_{\ol p_0})x^{i\si_{p_0}}.
\hspace{0.9cm}(4.152)\end{eqnarray*}
Comparing the coefficients of $x^{-\si_{p_0}}$, we obtain that
$$
4+\vt(\a',0)=b(2-\a_{p_0})-2(\a_{p_0}b+\vt(\a',0)),
\eqno(4.153)$$
holds for all $\a\in\G$ with $\a_0=0,\a_{p_0}\ne\a_{\ol p_0}$. This
is impossible. Thus $J'_0=\{0\}$. This proves Claim 10.
\par
By Claim 10, we also have $J_0=\{0\}$ and $p_0\in I_1$. Thus
by the sufficiency, there exists isomorphism $\wt\th:\AA'\cong
\wt\AA$ which extends $\wt\th:{\cal E}'\cong\wt{\cal E}$. Then
$\wt\th\cdot\th$ is an isomorphism $\AA\cong\wt\AA$ such that the
corresponding $\nu$ is the identity map, and so by Case 1(i), there exists
$g'\in{\cal G}'$ such that $g'_{\rm id}$ is a group isomorphism
$\G\cong\wt\G$. Thus $g'_\nu:\G\cong\G'$ is a group isomorphism as required.
This proves the theorem in this case.
\par
{\it Case 2}: $\G_0=\G'_0=\{0\}$.
First assume that $\iota_3>0$. Using Lemma 3.5(2), we have (4.65) and
as in Claim 2, there is a bijection $\nu:I_{1,3}\rar I'_{1,3}$.
By (4.80) (where now $\a_0=\b_0=0$), we can suppose as in (4.84),
$$
\tau(\a+\b+\si_1)=\tau(\a)+\tau(\b)+\si_{\nu(1)}+
\sum_{1\ne p\in I_{1,3}}k^{(p)}_{\a,\b}(\si_{\nu(p)}-\si_{\nu(1)}).
\eqno(4.154)$$
Then as proved before, we have $k^{(p)}_{\a,\b}=0$, and so
$$
\tau':\a\mapsto\tau(\a-\si_1)+\si_{\nu(1)},
\eqno(4.155)$$
can be uniquely extended to a group isomorphism $\tau':\G\rar\G'$ such that
$\tau(\a)=\tau'(\a)+\tau'(\si_1)-\si_{\nu(1)}$. To see how the proof is similar
as before, we assume that $c_\a=1$ in (4.65), $\nu={\rm Id}_{I_{1,3}}$
(the general case can be obtain by multiplying $\th$ by an element $g'_\nu\in
{\cal G}$). For convenience,
take $\a,\b\in\G\cap(\mbb{F}1_{[1]}+\mbb{F}1_{[\ol1]})$. We have
$$
[x^\a,x^\b]=(\a_1\b_{\ol1}-\a_{\ol1}\b_1)x^{\a+\b+\si_1}.
\eqno(4.156)$$
Thus
\begin{eqnarray*}& &
(((\tau'(\a))_1+(\tau'(\si_1))_1+1)((\tau'(\b))_{\ol1}+(\tau'(\si_1))_{\ol1}+1)\\ & &
-((\tau'(\a))_{\ol1}+(\tau'(\si_1))_{\ol1}+1)((\tau'(\b))_1+(\tau'(\si_1))_1+1))
x^{\tau'(\a)+\tau'(\b)+2\tau'(\si_1)-\si_1}\\ &=&[x^{\tau'(\a)+\tau'(\si_1)-\si_1},x^{\tau'(\b)+\tau'(\si_1)-\si_1}]\\ &=&(\a_1\b_{\ol1}-\a_{\ol1}\b_1)x^{\tau'(\a)+\tau'(\b)+2\tau'(\si_1)-\si_1}.
\hspace{7.5cm}(4.157)\end{eqnarray*}
Fixing $\a$ and comparing the coefficients, since $\tau'$ is a group
isomorphism, $(\tau'(\b))_1,(\tau'(\b))_{\ol1}$ must linearly depend
on $\b_1,\b_{\ol1}$. {}From this we obtain $(\tau'(\si_1))_1+1=0=
(\tau'(\si_1))_{\ol1}+1$. In general we can obtain $\tau'(\si_1)=\si_1$.
Hence $\tau'=\tau$. The rest of the proof is exactly the same as before except
that we do not have the second equation in (4.90) since $\a_0=\b_0=0$ in (4.97). Thus we have the
first case of (4.7). This proves the theorem in this case.
\par
Next assume that $\iota_3=0$. If $\ell_4+\ell_5=0$, there is nothing
to be proved (cf.~(4.120)). Assume that $\ell_4+\ell_5>0$.
Then by Lemma 3.5(3), we have (4.65) for
all $\a\in\G$. The rest of the proof is again similar as before.
\par
{\it Case 3}. $\G_0\ne\{0\}$ and $\G'_0=\{0\}$.
\par
Using Lemma 3.5(1) and (2), as the proof in Case 1,
there exists a bijection $\tau:\G\bs\G_{0,3}\to\G'\bs\G'_{0,3}$
such that $\th(x^\a)\in\mbb{F}x^{\tau(\a)}$ for $\a\in\G\bs\G_{0,3}$,
and there exists bijection $\nu:\wh I_{1,3}\to I'_{1,3}$ such that
$\th(x^{-\si_p})\in\mbb{F}x^{-\si_{\nu(p)}}$ for $p\in\wh I_{1,3}$.
Then as proved before, we see that $(\tau(\a))_{\nu(0)}$ and
$(\tau(\a))_{\ol{\nu(0)}}$ only linearly depend on $\a_0$.
This is impossible by the second condition of (2.15). Thus this case
does not occur.
$\qquad\Box$

\vspace{1cm}
\noindent{\Large \bf References}
%\hspace{0.5cm}

\begin{description}

\item[{[A]}] V. I. Arnold, {\it Mathematical Methods of Classical Mechanics}, second edition, Springer-Verlag New York, 1989.

\item[{[D]}] B, Dubrivin, ``Geometry of 2D topological field theories'' in  Integrable Systems and Quantum Groups (Montecatini Tenne, 1993), {\it Lecture Notes in Math.} {\bf 1020}, Springer-Verlag, Berlin, 1996, 120-384.

\item[{[HM]}] C. Hertling and Yu. Manin, Weak Frobenius manifolds, {\it Intern. Math. Res. Notice} (1999), no. 6, 277-280.

\item[{[K1]}] V. G. Kac, A description of filtered Lie algebras whose associated graded Lie algebras are of Cartan types, {\it Math. of USSR-Izvestijia} {\bf 8} (1974), 801-835.

\item[{[K2]}] V. G. Kac, Lie superalgebras, {\it Adv. Math.} {\bf 26} (1977), 8-96.

\item[{[K3]}] V. G. Kac, {\it Vertex algebras for beginners}, University lectures series, Vol {\bf 10}, AMS. Providence RI, 1996.

\item[{[K4]}] V. G. Kac, Classification of infinite-dimensional simple linearly compact Lie superalgebras, {\it Adv. Math.} {\bf 139} (1998), 1-55.

\item[{[O]}] J. Marshall Osborn, New simple infinite-dimensional Lie algebras of characteristic 0, {\it J. Algebra} {\bf 185} (1996), 820-835.

\item[{[OZ]}] J. Marshall Osborn and K. Zhao, Generalized Cartan type K Lie
algebras in characteristic 0, {\it Commun. Algebra} {\bf 25} (1997),
3325-3360.

\item[{[SF]}] H. Strade and R. Farnsteiner, {\it Modular Lie Algebras and Their Representations,} Marcel Derkker, Inc., 1988.

\item[{[SXZ]}] Y. Su, X. Xu and H. Zhang, Derivation-simple algebras and the structures of  Lie algebras of Witt type, {\it J. Algebra} {\bf 233} (2000), 642-662.

\item[{[X1]}] X. Xu, New generalized simple Lie algebras of Cartan type over a field with characteristic 0,  {\it J. Algebra} {\bf 224} (2000), 23-58.

\item[{[X2]}] X. Xu, Quadratic conformal superalgebras, {\it J. Algebra} {\bf 231} (2000), 1-38.

\end{description}

\end{document}